\documentclass[a4paper,11pt]{article}
\usepackage{latexsym,amssymb,enumerate,amsmath,epsfig,amsthm,dsfont,bm}
\usepackage[margin=1in]{geometry}
\usepackage{setspace,color}
\usepackage{tikz}
\usepackage{floatrow}
\usepackage{graphicx,subfigure}
\usepackage[ruled]{algorithm2e}
\usepackage{epstopdf}
\usepackage[multidot]{grffile}
\usepackage{comment}
\usepackage{multirow}
\newcommand{\D}{\mathcal{D}}
\newcommand{\x}{\mathbf{x}}
\newcommand{\y}{\mathbf{y}}
\newcommand{\p}{\mathbf{p}}

\newcommand{\R}{\mathbb{R}}
\newcommand{\rd}{\,\mathrm{d}}
\newcommand{\bflambda}{\bm{\lambda}}
\newcommand{\maF}{\mathcal{F}}

\makeatletter
\newcommand\figcaption{\def\@captype{figure}\caption}
\newcommand\tabcaption{\def\@captype{table}\caption}
\makeatother

\theoremstyle{plain}

\theoremstyle{definition}

\newtheorem{?}[Th]{Problem}

\begin{document}

\title{Fast Algorithms for Surface Reconstruction from Point Cloud}
\date{\vspace{-5ex}}
\author{
Yuchen He\footnote{royarthur@gatech.edu, School of Mathematics, Georgia Institute of Mathematics, Atlanta, USA},
Martin Huska \footnote{martin.huska@unibo.it, Department of Mathematics, University of Bologna, Italy},
Sung Ha Kang\footnote{kang@math.gatech.edu,School of Mathematics, Georgia Institute of Mathematics, Atlanta, USA.  Research is supported in part by Simons Foundation grant 584960.}, and
Hao Liu\footnote{hao.liu@math.gatech.edu, School of Mathematics, Georgia Institute of Mathematics, Atlanta, USA}
}
\maketitle

\begin{abstract}
We consider constructing a surface from a given set of point cloud data. We explore two fast algorithms to minimize the weighted minimum surface energy in [Zhao, Osher, Merriman and Kang, Comp.Vision and Image Under., 80(3):295-319, 2000].   An approach using Semi-Implicit Method (SIM) improves the computational efficiency through relaxation on the time-step constraint. An approach based on Augmented Lagrangian Method (ALM) reduces the run-time via an Alternating Direction Method of Multipliers-type algorithm, where each sub-problem is solved efficiently.  We analyze the effects of the parameters on the level-set evolution and explore the connection between these two approaches.  We present numerical examples to validate our algorithms in terms of their accuracy and efficiency.
\end{abstract}

\section{Introduction}

\begin{figure}[ht]
\centering
\begin{tabular}{ccc}
	(a)&(b)&(c)\\
	\includegraphics[width=1.3in,height=1.3in]{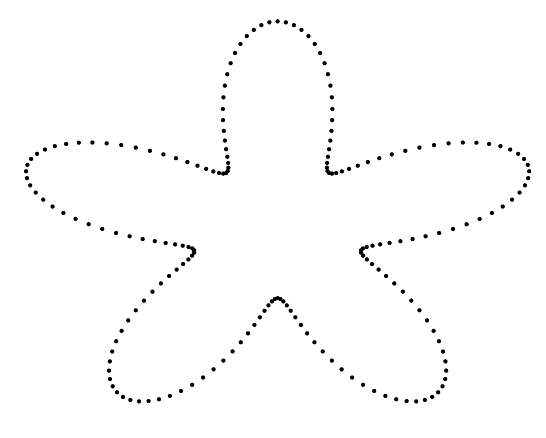}&
	\includegraphics[width=1.3in,height=1.3in]{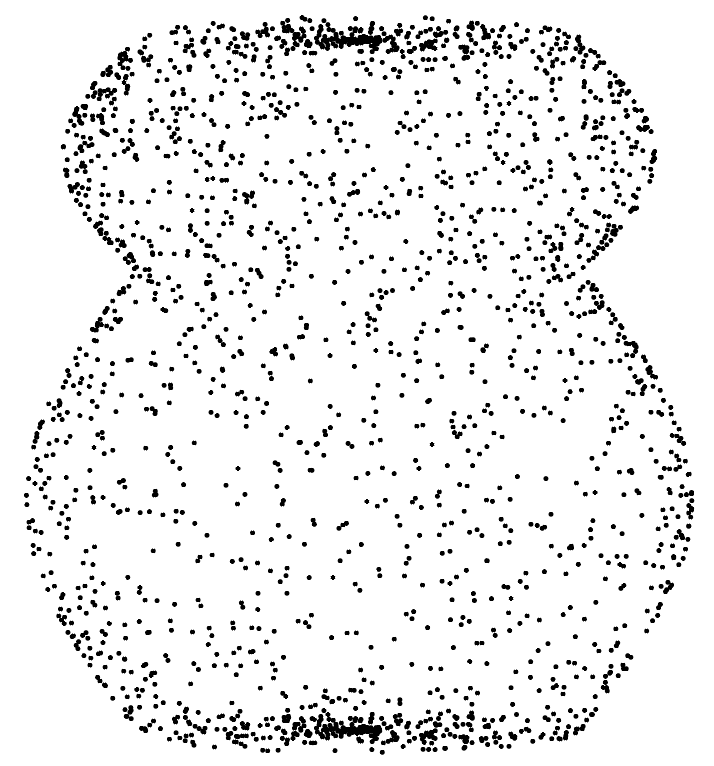}&
	\includegraphics[height=1.3in]{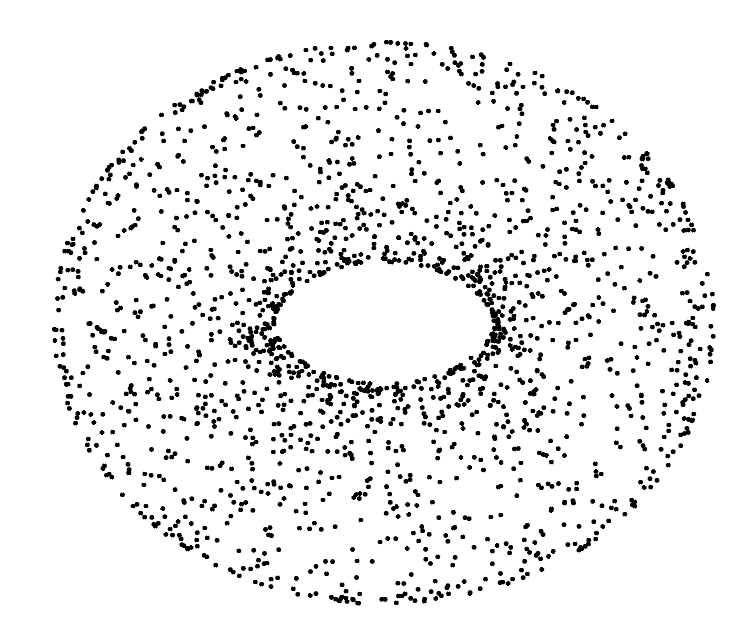}
\end{tabular}
\caption{Test point clouds. (a) Five-fold circle (200 points). (b) Jar (2100 points). (c) Torus (2000 points).}\label{fig::pcd}
\end{figure}

Acquisition, creation and processing of 3D digital objects is an important topic in various fields, e.g., medical imaging~\cite{khan2018single}, computer graphics~\cite{CLMM06,CT11}, industry \cite{bi2010advances}, and preservation of cultural heritage \cite{gomes20143d}.   A fundamental step is to reconstruct a surface  from a 3D scanned point cloud data~\cite{alexa2001point},  denoted as $\D\subseteq\mathbb{R}^m$ for $m=2$ or $3$, such as in Figure~\ref{fig::pcd}.

We focus on reconstructing an $m$-dimensional manifold $\Gamma$, a curve in $\mathbb{R}^2$  or a surface in $\mathbb{R}^3$, from the point cloud $\D$.   We assume only the point locations are given, and no other  geometrical information such as normal vectors at each point are known.  We explore fast algorithms for minimizing the following energy proposed in \cite{ZOM*98}:
\begin{align}
 E_p(\Gamma)=\left(\int_{\Gamma} |d(\x)|^p \rd\x\right)^{\frac{1}{p}},
  \label{eq.energy0}
\end{align}
where $d(\x)=\min_{\y\in \D}\left\{|\x-\y|\right\}$ is the distance from any point $\x\in\R^m$ to $\D$, $p $ is a positive integer, and $\rd\x$ is the surface area element.  This energy is the $p$-norm of the distance function restricted in $\Gamma$.
In \cite{ZOM*98}, the authors used the fast sweeping scheme to solve Eikonal equation, and the Euler-Lagrange equations are solved by gradient descent algorithm.

Among many ways to represent the  underlying surface, e.g., (moving) least square projection \cite{ABCO*03,OGG09}, radial basis function \cite{DTS01,CBC*01,CBC*03,CLMM06,HLZ*09}, poisson reconstruction \cite{KBH06,BKBH09,LWCC10,KH13,ESTS15}, we use the level set method as in \cite{EZL*12,ZOF01, ZOM*98}.   The level set formulation allows topological changes and self-intersection during the evolution~\cite{OS88,OF01}, and gain popularity in many applications \cite{chan2001active,tsai2003shape,sethian1999level} just to mention a few.
We represent the surface as a zero level set of $\phi:\mathbb{R}^{m+1}\to\mathbb{R}$:
\begin{align*}
\Gamma=\phi^{-1}(0)=\{\x\in\mathbb{R}^m\mid \phi(\x)=0\}.
\end{align*}	
%and we search for a $\phi$ minimizing the following energy:
%\begin{eqnarray}
%  E_p(\phi)=\left(\int_{\{\phi=0\}} |d(\x)|^p \d\x\right)^{\frac{1}{p}}\;
%  \label{eq.energy0_5}
%\end{eqnarray}

There are various related works on surface reconstruction from point cloud data: a convection model proposed in \cite{ZOF01}, a data-driven logarithmic prior for noisy data in \cite{SK04}, using surface tension to enrich the Euler-Lagrange equations in \cite{HM16}, and using principal component analysis to reconstruct curves embedded in sub-manifolds in \cite{liu2017level}.   A semi-implicit scheme is introduced in \cite{smereka2003semi} to simulate the curvature and surface diffusion motion of the interface. In  \cite{LPZ13}, the authors defined the surface via a collection of anisotropic Gaussians centered at each entry of the input point cloud, and used TVG-L1 model  \cite{bresson2007fast} for minimization.  A similar strategy addresses an $\ell_0$ gradient regularization model proposed in \cite{LLY*18}.
Some models incorporate additional information. In \cite{LWQ09}, the authors proposed a novel variational model, consisting of the distance, the normal, and the smoothness term.  Euler's Elastica model %\cite{TD17,BTZ17}
is incorporated for surface reconstruction in \cite{SWT*12} where graph cuts algorithm is used.   The model in \cite{EZL*12} extends the active contours segmentation model to 3D and implicitly allows controlling the curvature of the level set function.

In this paper, we explore fast algorithms to minimize the weighted minimum surface energy (\ref{eq.energy0}) for $p=1$ and 2.   We propose a Semi-Implicit Method (SIM) to relax the time-step constraint for $p=2$, and an Augmented Lagrangian Method (ALM) based on the alternating direction method of multipliers (ADMM) approach for $p=1$.  These algorithms minimize the weighted minimal surface energy~(\ref{eq.energy0}) with high accuracy and superior efficiency.   We analyze the behavior of ALM in terms of parameter choices and explore its connection to SIM.   Various numerical experiments are presented to discuss the effects of the algorithms.

We organize this paper as follows.  In Section~\ref{SEC:model}, we present the two methods: SIM and ALM, and explore their connection.
Numerical experiments are presented in Section \ref{SEC:exp}, and effects of the parameters are discussed in Subsection~\ref{SEC::ALMbehavior}.   We conclude our paper in Section \ref{SEC:concl}.

\section{Proposed Algorithms}
\label{SEC:model}

Let $\Omega\subset\mathbb{R}^{m}$ ($m=2$ or $3$) denote a bounded domain containing the given $m$-dimensional point cloud data, $\D$.  Using the level-set formulation for $\Gamma$,  the $d$-weighted minimum surface energy (\ref{eq.energy0}) can be rewritten as:
\begin{eqnarray}
  E_p(\phi)=\left(\int_{\Omega} |d(\x)|^p \delta(\phi)|\nabla\phi|\rd\x\right)^{\frac{1}{p}}.
  \label{eq.energy1}
\end{eqnarray}
Here $\delta(x)$ is the Dirac delta function which takes $+\infty$ when $x=0$, and $0$ elsewhere.  Compared to (\ref{eq.energy0}), this integral is defined on $\Omega$, which makes the computation flexible and free from explicitly tracking $\Gamma$.
We use $p=2$ for SIM introduced in Section~\ref{SUBS:L2_scheme}, and $p=1$ for ALM in Section~\ref{SUBS:L1_scheme}.  In general, $p=2$ is a natural choice, since it provides better stability and efficiency for a semi-implicit type PDE-based method.  For ALM, we explore $p=1$ to take advantage of an aspect of fast algorithm in ADMM setting such as shrinkage,  similarly to the case in \cite{BTZ17}.   Visually, the numerical results of surface reconstruction are similar for $p=1$ or $p=2$ (see Section \ref{SEC:exp}).

\subsection{Semi-Implicit Method (SIM) to minimize $\mathbf{E_2}$}
\label{SUBS:L2_scheme}

We introduce a gradient-flow-based semi-implicit method to Minimize %(\ref{eq.energy1}) for $p=2$, i.e.,
\begin{equation}
E_2(\phi)=\left(\int_{\Omega} d(\x)^2 \delta(\phi)|\nabla\phi| \rd\x\right)^{\frac{1}{2}}.
\label{eq.split}
\end{equation}
Following~\cite{ZOM*98}, the first variation of $E_2(\phi)$ with respect to $\phi$ is characterized as a functional:
\begin{eqnarray*}
\left\langle \frac{\partial E_2(\phi)}{\partial\phi},v \right\rangle &= &
\int_{\Omega}\frac{1}{2}\delta(\phi)
\left[ \int_{\Omega} d^2(\x)\delta(\phi)|\nabla\phi|\rd\x\right]^{-1/2}
\nabla\cdot \left[ d^2(\x)\frac{\nabla\phi}{|\nabla \phi|}\right]v \rd\x
\label{eq.gradient}
\end{eqnarray*}
for any test function $v$ from the Sobolev space $H^1$.
Minimizing (\ref{eq.split}) is equivalent to finding the critical point $\phi$ such that
$%\begin{equation}$
\left\langle \frac{\partial E_2(\phi)}{\partial\phi},v\right\rangle=0,
\forall v\in H^1. $
%\label{eq.euler}
%\end{equation}
This is associated with solving the following initial value problem:
\begin{equation}
\begin{cases}
\displaystyle{\frac{\partial \phi}{\partial t}=\bar{f}(d,\phi) \nabla\cdot \left[ d^2(\x)\frac{\nabla\phi}{|\nabla \phi|}\right]} ,\\
\phi(\x,0)=\phi^0,
\end{cases}
\label{eq.ivp}
\end{equation}
where $\phi^0$ is an initial guess for the unknown $\phi$, and
$
\displaystyle{\bar{f}(d,\phi)=\frac{1}{2}\delta(\phi)\left[ \int_{\Omega} d^2(\x)\delta(\phi)|\nabla\phi|\rd\x\right]^{-1/2}.}
$ The steady state solution of (\ref{eq.ivp}) gives a minimizer $\phi^*$ of $E_2(\phi)$.

Here the delta function $\delta$ is realized as the derivative of the one dimensional
	Heaviside function $H:\mathbb{R}\to\{0,1\}$.
%taking value $1$ on positive numbers, $-1$ on negative ones, and $1/2$ at $0$.
We adopt the smooth approximation of $H(\phi)$ as in \cite{Bra00}:
\begin{equation}
	H(\phi) \approx H_\varepsilon(\phi) = \frac{1}{2}+\arctan(\phi/\varepsilon)/\pi \;\; \text{ and } \;\;
	\delta(\phi)\approx H^\prime_\varepsilon(\phi) = \frac{\varepsilon}{\pi(\varepsilon^2+\phi^2)} \;
\label{eq.H_eps}
\end{equation}
with  $\varepsilon>0$ as the smoothness parameter.
Then $\bar{f}$ is approximated  by its smoothed version $f$ expressed as
\begin{equation*}
f(d,\phi)=\frac{1}{2}\frac{\varepsilon}{\pi(\varepsilon^2+\phi^2)}\left[ \int_{\Omega} d^2(\x)\frac{\varepsilon}{\pi(\varepsilon^2+\phi^2)}|\nabla\phi|\rd\x\right]^{-1/2}.
\end{equation*}

We add a stabilizing diffusive term $-\beta\Delta\phi$
	for $\beta>0$ on both sides of the PDE in~(\ref{eq.ivp}) to consolidate the computation, similarly to \cite{smereka2003semi}:
\begin{equation}
\frac{\partial \phi}{\partial t} - \beta\Delta \phi=-\beta\Delta \phi +f(d,\phi) \nabla\cdot \left[ d^2(\x)\frac{\nabla\phi}{|\nabla \phi|}\right].
\label{eq.ivp2}
\end{equation}
Employing a semi-implicit scheme, we solve $\phi$ from (\ref{eq.ivp2}) by iteratively updating $\phi^{n+1}$ using $\phi^{n}$ via the following equation:
\begin{equation}
\frac{\phi^{n+1}}{\Delta t}-\beta\Delta\phi^{n+1}=\frac{\phi^n}{\Delta t}-\beta\Delta\phi^n+f(d,\phi^n) \nabla\cdot \left[ d^2(\x)\frac{\nabla\phi^n}{|\nabla \phi^n|}\right],
\label{eq.phi.semi}
\end{equation}
where $\Delta t$ is the time-step.    This equation can be efficiently solved by the Fast Fourier Transform (FFT).  Denoting the discrete Fourier transform by $\maF$ and its inverse by $\maF^{-1}$,  we have:
\begin{eqnarray*}
	&\maF (\phi)(i\pm1,j)=e^{\pm2\pi\sqrt{-1}(i-1)/M}\maF (\phi)(i,j), &\maF (\phi)(i,j\pm 1)=e^{\pm2\pi\sqrt{-1}(j-1)/N}\maF (\phi)(i,j).
\end{eqnarray*}
Accordingly, the discrete Fourier transform of $\Delta\phi$ is
\begin{equation*}
\maF (\Delta \phi)(i,j)=\left[2\cos(\pi\sqrt{-1}(i-1)/M)+2\cos(\pi\sqrt{-1}(j-1)/N)-4\right]\maF \phi(i,j) .
\end{equation*}
Here the coefficient in front of $\maF \phi(i,j)$ represents the diagonalized discrete Laplacian operator in the frequency domain.
Let $g_1$ be the right side of (\ref{eq.phi.semi}), then the solution $\phi^{n+1}(i,j)$ of (\ref{eq.phi.semi}) is computed via
\begin{equation}
  \phi^{n+1}(i,j)=\maF^{-1}\left(\frac{\maF(g_1)(i,j)}{ \left(1-\beta\Delta t\left[2\cos(\pi\sqrt{-1}(i-1)/N)+2\cos(\pi\sqrt{-1}(j-1)/N)-4\right]\right)}\right).
\label{eq.FFT_SIM}
\end{equation}

As for the stopping criterion, we exploit the mean relative change of the weighted minimum surface energy~(\ref{eq.energy0}). At the $n^{\text{th}}$ iteration, the algorithm terminates if
\begin{equation}\label{eq::errordef}
	\frac{|\bar{e}^k_{n-1}-\bar{e}^k_n|}{\bar{e}^k_n} < 10^{-4}, \;\;
	\text{ where } \;\; \bar{e}_n^k=\frac{1}{k}\sum_{i=n-k}^n E_p(\phi^i).
%\label{eq.stop}
\end{equation}
Here the quantity $\bar{e}^k_n$ represents the average of the energy values computed from the $(n-k)^{\text{th}}$ to the $n^{\text{th}}$ iteration for some $k\in\mathbb{N}$, $k\geq 1$.  We fix $k=10$, and set $p=2$ for SIM.  We summarize the main steps of SIM in Algorithm~\ref{alg:SIM}.

\begin{algorithm}
	\SetKwInOut{KwIni}{Initialization}
	\KwIni{$d$, $\phi^0$ and $n=0$.}

	\While {the stopping criterion \eqref{eq::errordef} with $p=2$ is greater than $10^{-4}$}{
		Update $\phi^{n+1}$ from $\phi^{n}$ solving (\ref{eq.FFT_SIM})\;
%		Reinitialize $\phi^{n+1}$ \;
		Update $n\gets n+1$\;}
	\KwOut{$\phi^{n}$ such that $\{\phi^n=0\}$ approximates $\{\phi^*=0\}$.}
	\caption{SIM for the weighted minimum surface (\ref{eq.split})} \label{alg:SIM}
\end{algorithm}

\subsection{Augmented Lagrangian Method (ALM) to Minimize $\mathbf{E_1}$}
\label{SUBS:L1_scheme}

In this section, we present an augmented Lagrangian-based method to minimize the weighted minimum surface energy~\eqref{eq.energy1} for $p=1$, i.e.,
\begin{equation}
E_1(\phi)=\int_{\Omega}d(\x) \delta(\phi)|\nabla\phi| \rd\x.	
\label{eq.energy2}
\end{equation}

For the non-differentiable term $|\nabla\phi|$ in \eqref{eq.energy2},
we utilize the variable-splitting  and  introduce an auxiliary variable $\p=\nabla\phi$.
We rephrase the minimization of $E_1(\phi)$ as a constrained optimization problem:
\begin{align}
\{\phi^*,\p^*\} \gets \arg\min\limits_{\phi,\p}
\int_{\Omega}\frac{\varepsilon d |\p|}{\pi(\varepsilon^2+\phi^2)}\rd\x,\;\;\;
 \mathrm{subject\;to} \;\; \p=\nabla\phi,
\label{eq.energy2_constr}
\end{align}
here we replace $\delta(\phi)$ by its smooth approximation $H^\prime_\varepsilon(\phi)$ as in \eqref{eq.H_eps}.
To solve problem \eqref{eq.energy2_constr}, we formulate the
augmented Lagrangian function:
\begin{align}
\mathcal{L}(\phi,\p,\bflambda; r)=
\int_{\Omega}\frac{\varepsilon d |\p|}{\pi(\varepsilon^2+\phi^2)} \rd\x
+ \frac{r}{2}\int_{\Omega}|\p-\nabla\phi|^2 \rd\x
+ \int_{\Omega}\bflambda\cdot(\p-\nabla\phi) \rd\x,
\label{eq.lag}
\end{align}
where $r>0$ is a scalar penalty parameter and $\bflambda: \R^{m+1}\to\R^{m+1}$ represents
the Lagrangian multiplier.
Minimizing~(\ref{eq.lag}) amounts to considering the following saddle-point problem:
\begin{align}
\mathrm{Find}&\; (\phi^*,\p^*,\bflambda^*)\in\R\times\R^{m+1}\times\R^{m+1} \nonumber\\
\mathrm{s.t.}&\; \mathcal{L}(\phi^*,\p^*,\bflambda; r)\le\mathcal{L}(\phi^*,\p^*,\bflambda^*; r)
\le\mathcal{L}(\phi,\p,\bflambda^*; r); \nonumber \\
&\;\forall (\phi,\p,\bflambda)\in\R\times\R^{m+1}\times\R^{m+1}.
\label{eq.saddle_point}
\end{align}
Given $\phi^n$, $\p^n$, and $\bflambda^n$, for $n=0,1,2,\cdots$,
the $(n+1)^{\text{th}}$ iteration of an ADMM-type
	algorithm for \eqref{eq.saddle_point} consists of solving a series of sub-problems:
\begin{align}
\phi^{n+1}	&=\arg\min_{\phi}\mathcal{L}(\phi,\p^n,\bflambda^n;r);\label{eq.sub1}\\
\p^{n+1}	&=\arg\min_{\p}\mathcal{L}(\phi^{n+1},\p,\bflambda^n;r);\label{eq.sub2}\\
\bflambda^{n+1}&=\bflambda^{n}+r\left(\p^{n+1}-\nabla\phi^{n+1}\right).\label{eq.sub3}
\end{align}

Each sub-problem can be solved efficiently.
First, we find the minimizer of the sub-problem (\ref{eq.sub1}) by solving its Euler-Lagrange equation:
\begin{equation}
-r\Delta\phi^{n+1}=\frac{2d\varepsilon|\p^n|\phi^n}{\pi(\varepsilon^2+(\phi^n)^2)^2}-\nabla\cdot(r\p^n+\bflambda^n).
\label{eq.sub1_PDE}	
\end{equation}
Following~\cite{BTZ17}, we introduce a frozen-coefficient term $\eta\phi$, for $\eta>0$,
	on both sides of (\ref{eq.sub1_PDE}) to stabilize the computation; thus, (\ref{eq.sub1}) is solved using the following equation:
\begin{equation}
\eta\phi^{n+1}-r\Delta\phi^{n+1}=\eta\phi^{n}+\frac{2d\varepsilon|\p^n|\phi^n}{\pi(\varepsilon^2+(\phi^n)^2)^2}-\nabla\cdot(r\p^n+\bflambda^n).
\label{eq.sub1_PDE1}	
\end{equation}
Here $\Delta$ is the Laplacian operator, and we solve this via FFT,  similar to  (\ref{eq.FFT_SIM}) for SIM.
Thus, the
$\phi$ sub-problem is solved via:
\begin{equation}
  \phi^{n+1}(i,j)=\maF^{-1}\left(\frac{\maF(g_2)(i,j)}{ \left(\eta-r\left[2\cos(\pi\sqrt{-1}(i-1)/N)+2\cos(\pi\sqrt{-1}(j-1)/N)-4\right]\right)}\right).
\label{eq.FFT_ALM}
\end{equation}

Second, the $\p$ sub-problem \eqref{eq.sub2} is equivalent to a weighted Total Variation (TV) minimization,
	whose solution admits a closed-form expression using the shrinkage operator~\cite{tai2011fast}. Explicitly, the updated $\p^{n+1}$  is computed via:
\begin{align}
\p^{n+1} = \max\left\{0,1-\frac{d\,\varepsilon}{\pi(\varepsilon^2+(\phi^{n+1})^2)|r\nabla\phi^{n+1}-\bflambda^n|}\right\}
\left(\nabla\phi^{n+1}-\frac{\bflambda^n}{r}\right).
\label{eq.sub2_sol}
\end{align}

Finally, the Lagrangian multiplier $\bflambda$ is updated by~(\ref{eq.sub3}). The stopping criterion for the ALM iteration is the same as that for SIM (\ref{eq::errordef}), but with $p=1$.   We summarize the main steps of ALM in Algorithm~\ref{alg:ALM}.

\begin{algorithm}
	\SetKwInOut{KwIni}{Initialization}
	\KwIni{$d$, $\phi^0$, $\mathbf{p}^0$, $\bflambda^0$, and $n=0$.}
	\While {the stopping criterion \eqref{eq::errordef} with $p=1$ is greater than $10^{-4}$}{
		Update $\phi^{n+1}=\arg\min_{\phi}\mathcal{L}(\phi,\mathbf{p}^{n},\bflambda^n;r)$ via  (\ref{eq.FFT_ALM}) \;
%		Reinitialize $\phi^{n+1}$ \;
		Update $\mathbf{p}^{n+1}=\arg\min_{\mathbf{p}}\mathcal{L}(\phi^{n+1},\mathbf{p},\bflambda^n;r)$ via (\ref{eq.sub2_sol})\;
		Update $\bflambda^{n+1}=\bflambda^{n}+r(\mathbf{p}^{n+1}-\nabla\phi^{n+1})$\;
		Update $n\gets n+1$\;}
	\KwOut{$\phi^{n}$ such that $\{\phi^n=0\}$ approximates $\{\phi^*=0\}$.}
	\caption{ALM for the weighted minimum surface (\ref{eq.energy2})}
	\label{alg:ALM}
\end{algorithm}

%%%%%
\subsection{Connection between SIM and ALM algorithms}\label{SEC::conn}

Note that both SIM and ALM involve solving  elliptic PDEs of the form:
\begin{equation}
  a\phi-b\Delta\phi=g,
  \label{eq.general}
\end{equation}
for some constants $a,b>0$, and a function $g$  defined on $\Omega$.
For SIM, it is equation~(\ref{eq.phi.semi}):
\[
\underbrace{\frac{1}{\Delta t}}_{a}\phi^{n+1}-\underbrace{\beta}_{b}\Delta\phi^{n+1}=\underbrace{\frac{\phi^n}{\Delta t}-\beta\Delta\phi^n+f(d,\phi^n) \nabla\cdot \left[ d^2(\x)\frac{\nabla\phi^n}{|\nabla \phi^n|}\right]}_{g},
\]
and for ALM, it is equation~(\ref{eq.sub1_PDE1}):
\[
\underbrace{\eta}_{a}\phi^{n+1}-\underbrace{r}_{b}\Delta\phi^{n+1}=\underbrace{\eta\phi^{n}+\frac{2d\varepsilon|\p^n|\phi^n}{\pi(\varepsilon^2+(\phi^n)^2)^2}-\nabla\cdot(r\p^n+\bflambda^n)}_{g}.
\]

\begin{figure}
\centering
\includegraphics[scale=0.25]{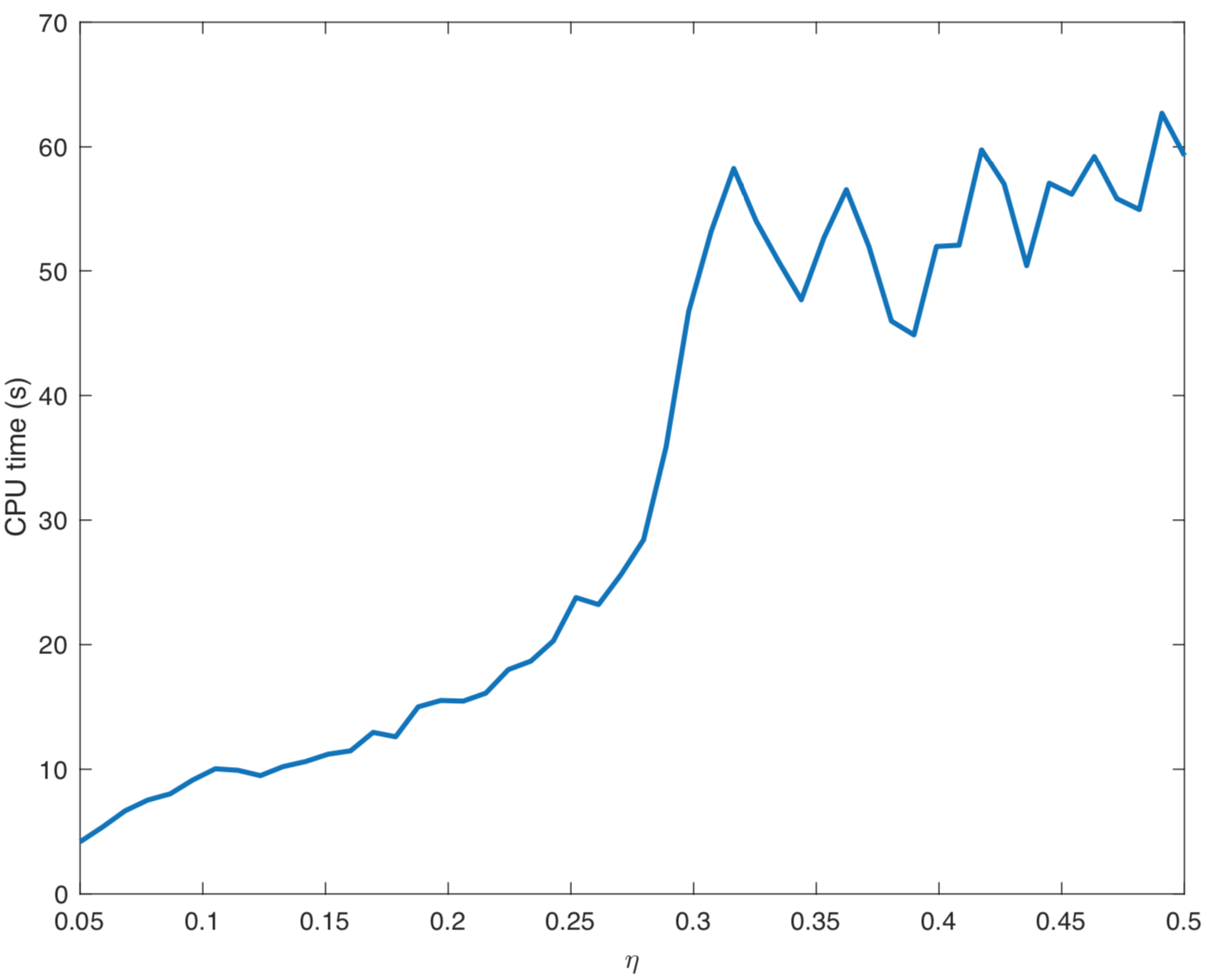}	
\caption{The CPU-time (s) of ALM until convergence, for the five-fold circle point cloud in Figure~\ref{fig::pcd}~(a). Here $r=\varepsilon=1$ and $\eta$ varies from $0.05$ to $0.5$.
The connection between SIM and ALM indicates that large $\eta$ slows down ALM.   In this  graph,
as $\eta$ increases, the time required to reach the convergence increases.}\label{fig::etacompute}
\end{figure}

We remark  interesting connections between SIM and ALM.
First, both methods have stabilizing terms but in different positions on the left side of (\ref{eq.general}). For SIM,  it is $-\beta \Delta\phi$, while for ALM, it is $\eta \phi$.
Second, relating  the coefficients of $\phi$,  $1/\Delta t$ in SIM gives insight to the effect of $\eta$ in ALM.  In general, a large $\eta$ slows down the convergence of ALM, while a small $\eta$ accelerates it (as the effect of $\frac{1}{\Delta t}$ on SIM).
Figure~\ref{fig::etacompute} shows convergence behaviors of ALM for different $\eta$, using the five-fold circle point cloud in Figure~\ref{fig::pcd}~(a).   It displays the CPU time~(in seconds) for $r=1$, $\varepsilon=1$, and $\eta$ varying from $0.05$ to $0.5$.
Note that as $\eta$ increases, the time required to reach the convergence increases almost quadratically at first, then stays around the same level.
Third, the correspondence between $b=\beta$ in SIM, and $b=r$ in ALM allows another interpretation of the parameter $r$ in ALM. In SIM, a large $\beta$ smears the solution and avoids discontinuities or sharp corners, and for ALM, large $r$ also allows to pass through fine details.  Figure~\ref{fig.fivecircle} in Section \ref{SEC:exp} presents more details, where we experiment with different $r$ and $\varepsilon$ values for the five-fold circle point cloud shown in Figure~\ref{fig::pcd}~(a).

\section{Numerical Implementations, Experiments and Effects of Parameters}
\label{SEC:exp}

In this section, we describe the implementation details and present numerical experiments.
For both SIM and ALM, we vary $\varepsilon$ from $0.5$ to $1$. For SIM, we use $\Delta t=500$.    When $\D$ is in 2D, we set $\beta=0.1$, and $\beta=0.01$ for 3D. For ALM, $\eta$ ranges from $0.05$ to $1$, and $r$ from $0.5$ to $2$.
%We discuss the interaction between $r$ and $\varepsilon$  in Section~\ref{SEC::ALMbehav} and the effect of $\eta$ in Section~\ref{SEC::conn} within wider ranges.

%To compare the CPU time with the original method \cite{ZOM*98}, we mimick it
%	by solving \eqref{eq.ivp2} in its fully explicit form.
	
The code is written in \textsc{Matlab} and executed without additional machine support,
	e.g. parallelization or GPU-enhanced computations.
All the experiments are performed on Intel\textsuperscript{\textregistered} Core\texttrademark
	4-Core 1.8GHz (4.0GHz with Turbo) machine, with 16 GB/RAM and Intel\textsuperscript{\textregistered} UHD Graphics 620 graphic card under Windows OS.
The contours and isosurfaces are displayed using \textsc{Matlab} visualization engine. No post-processing, e.g., smoothing or sharpening, is applied.

\subsection{Implementation Details}
\label{SEC:num}

We illustrate the details for planar point clouds, i.e., $\D\subseteq\R^2$, and the extension to $\mathbb{R}^3$ is straightforward.
Let the computational domain $\Omega = [0,M]\times[0,N]$, $M,N>0$, be discretized by a Cartesian grid with $\Delta x=\Delta y=1$.
For any function $u$ (or a vector field $\mathbf{v}=(v^1,v^2)$) defined on $\Omega$, we use  $u_{i,j}$ or $u(i,j)$ to denote $u(i\Delta x,i\Delta y)$.  We use the usual backward and forward finite difference schemes:
\begin{eqnarray*}
&&\partial_1^-u_{i,j}=\begin{cases}u_{i,j}-u_{i-1,j},~1<i\leq M;\\
u_{1,j}-u_{M,j},~i=1.	
\end{cases} \;\;\;
 \partial_1^+u_{i,j}=\begin{cases}
u_{i+1,j}-u_{i,j},~1\leq i<M-1;\\
u_{1,j}-u_{M,j},~i=M.	
\end{cases}
\\
&&\partial_2^-u_{i,j}=\begin{cases}u_{i,j}-u_{i,j-1},~1<j\leq N;\\
u_{i,1}-u_{i,N},~j=1.	
\end{cases} \;\;\;
 \partial_2^+u_{i,j}=\begin{cases}u_{i,j+1}-u_{i,j},1\leq j<N-1;\\
	u_{i,1}-u_{i,N},~j=N.
\end{cases}
\end{eqnarray*}
The gradient, divergence and the Laplacian operators are approximated as follows:
\begin{align*}
\nabla u_{i,j}&=((\partial_1^-u_{i,j}+\partial_1^+u_{i,j})/2, (\partial_2^-u_{i,j}+\partial_2^+u_{i,j})/2); \\%\label{eq.grad},\\
\nabla\cdot\mathbf{v}_{i,j} &= (\partial_1^+v_{i,j}^1 + \partial_1^-v_{i,j}^1)/2 + (\partial_2^+v_{i,j}^2 + \partial_2^-v_{i,j}^2)/2;\\
\Delta u_{i,j}&= \partial_1^+u_{i,j}-\partial_1^-u_{i,j}+ \partial_2^+u_{i,j}-\partial_2^-u_{i,j}.%\label{eq.laplace}
\end{align*}

The distance function $d$ is computed once at the beginning and no update is needed. It satisfies an Eikonal equation:
\begin{equation}
\begin{cases}
  |\nabla d|=1 \mbox{ in } \Omega,\\
  d(\x)=0 \mbox{ for } \x\in \D,
\end{cases}
\label{eq.eikonal}
\end{equation}
and discretizing (\ref{eq.eikonal}) via the Lax-Friedrich scheme leads to an updating formula:
\begin{eqnarray}
  &&d_{i,j}^{n+1}=\frac{1}{2}\Bigg(1-|\nabla d^n_{i,j}|+ \frac{d^n_{i+1,j}+d^n_{i-1,j}}{2} +  \frac{d^n_{i,j+1}+d^n_{i,j-1}}{2}\Bigg).
  \label{eq.updated}
\end{eqnarray}
 We solve (\ref{eq.updated}) using the fast sweeping method~\cite{kao2004lax} with complexity $O(G)$ for $G$ grid points.

Keeping $\phi^n$ to be a signed distance function during the iteration improves the stability of level-set-based algorithms.  We reinitialize $\phi^n$ at the $n^{\text{th}}$ iteration by  solving the following PDE:
\begin{equation}
\begin{cases}
\psi_{\tau}+\mbox{sign}(\psi)(|\nabla \psi|-1)=0,\\
\psi(\x,0) = \phi^n.
\end{cases}
\label{eq.reinitial}
\end{equation}
Here subscript $\tau$ represents artificial time and partial derivative, and $\text{sign}:\mathbb{R}\to\{-1,0,1\}$ is the sign function. In practice, $\psi^n$ being a signed distance function near $0$-level-set is important; thus, it is sufficient to solve (\ref{eq.reinitial}) for a few steps.  We fix  10 steps of reinitialization throughout this paper.

\subsection{Numerical Experiments of 2D and 3D Point Clouds}

\begin{figure}
\centering
\includegraphics[width=5in]{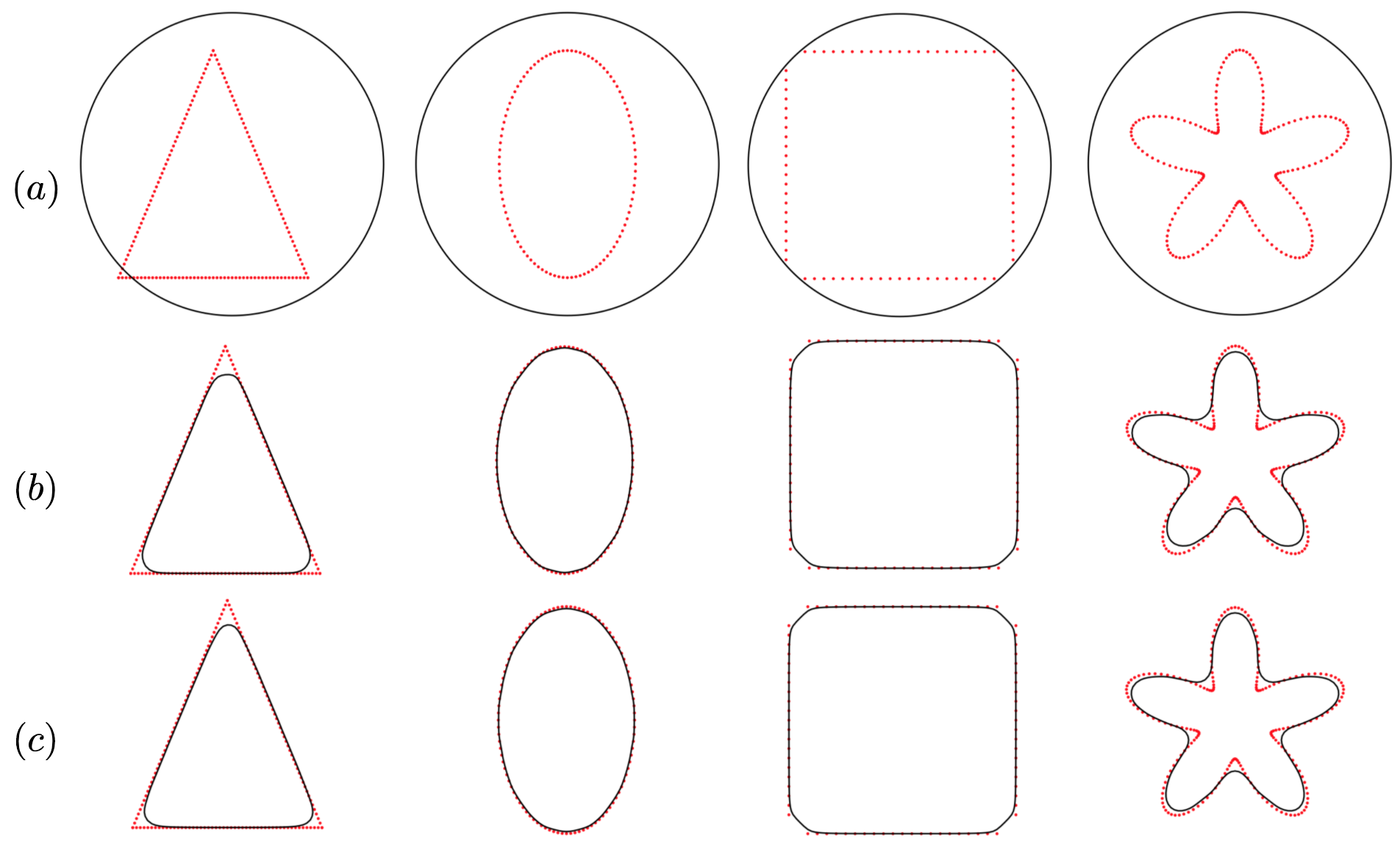}
\caption{The test point clouds: triangle with 150 number of points, ellipse with 100 points, square with 80 points, and five-fold-circle with 200 points.  (a) The top row, identical initial condition applied to SIM and ALM for different $\D$. (b) The middle row, the results obtained by SIM. (c) The bottom row, the results obtained by ALM using $r=1.5$. Both methods give compatible results. }\label{fig.2d}
\end{figure}

The first experiment, Figure~\ref{fig.2d}, is reconstruction of planar curves from 2D point clouds confined within a square $\Omega=[0,100]^2\subset\mathbb{R}^2$.  We generate the data using four different shapes: a triangle, an ellipse, a square whose corners are missing, and a five-fold circle.  For these cases, we use a centered circle with radius $30$ as the initial guess, shown in Figure~\ref{fig.2d} (a). Figure~\ref{fig.2d} (b) and (c) display the given $\D$, as well as the curves identified by SIM and ALM with $r=1.5$, respectively. Both methods produce comparably accurate results.  In the triangle example, corners get as close as the approximated delta function (with parameter $\varepsilon$) allows for both methods.  The ellipse and square results fit very closely to the given point clouds.   For the five-fold-circle, there is a slight difference in how the curve fits the edges, yet the results are very compatible.

Table \ref{tab.2d.cpu} shows the CPU time (in seconds) for SIM, ALM using $r=0.5, 1,1.5$, and $2$, and the explicit method in \cite{ZOM*98} with $\Delta t=20$ on the same data sets. With proper choices of $r$, ALM outperforms the other methods in terms of computational efficiency. SIM is stable without any dependency on the choice of parameters, and its run-times are comparable to the best performances of ALM in most cases. Both methods are faster than the explicit method in all the examples, and computational efficiency is similar between ALM and SIM.
\begin{table}[h]
	\begin{center}
		\begin{tabular}{|c|c|c|c|c|c|c|}
			\hline
			Object&ALM($r=0.5$)&ALM($r=1)$ &ALM($r=1.5$)&ALM($r=2$)& SIM & \cite{ZOM*98}\\\hline
			Triangle& $-$& $1.45$ & $1.31$ & $1.48$ &$1.50$ & $5.25$\\\hline
			Ellipse& $1.22$ & $1.03$& $1.33$ & $1.37$ &$1.49$ & $3.89$\\\hline
			Square& $-$& $-$& $0.94$ & 1.20&$1.09$ & $2.07$\\\hline	
			Five-fold circle &  $0.83$&$1.44$ & $1.86$& $1.22$&$1.96$ & $4.18$\\\hline
		\end{tabular}
	\end{center}
	\caption{CPU time (s) for SIM, ALM using $r=0.5, 1,1.5$, and $2$, and the explicit method in \cite{ZOM*98} with $\Delta t=20$ for the  point cloud data sets in Figure~\ref{fig.2d}.  Both SIM and ALM shows fast convergence. }
	\label{tab.2d.cpu}
\end{table}

\begin{figure}
  \centering
  \begin{tabular}{cc}
  (a)&(b)\\	
  \includegraphics[width=1.5in,height=1.5in]{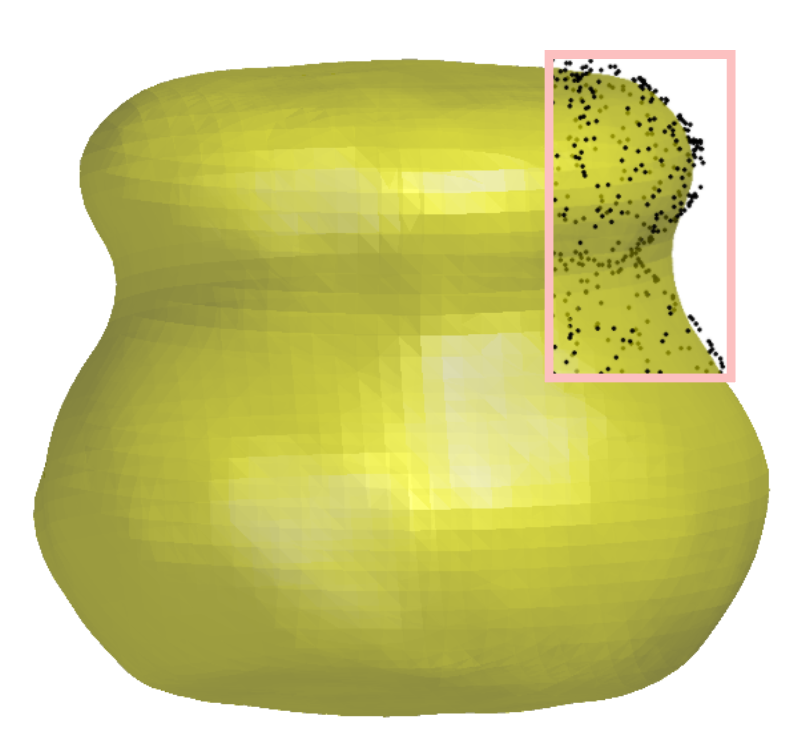}&
  \includegraphics[width=1.5in,height=1.5in]{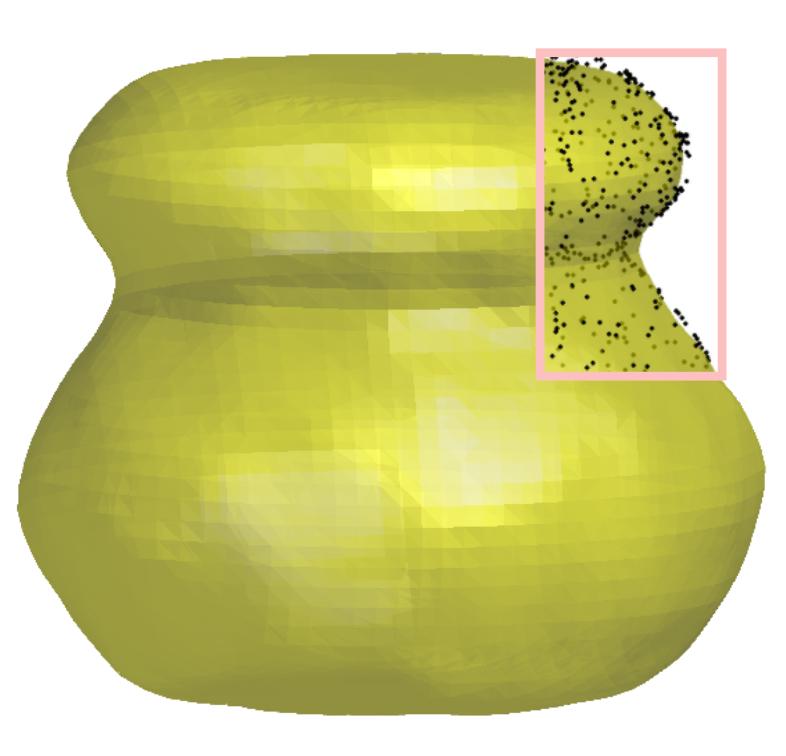}\\
  (c)&(d)\\
  \includegraphics[scale=0.3]{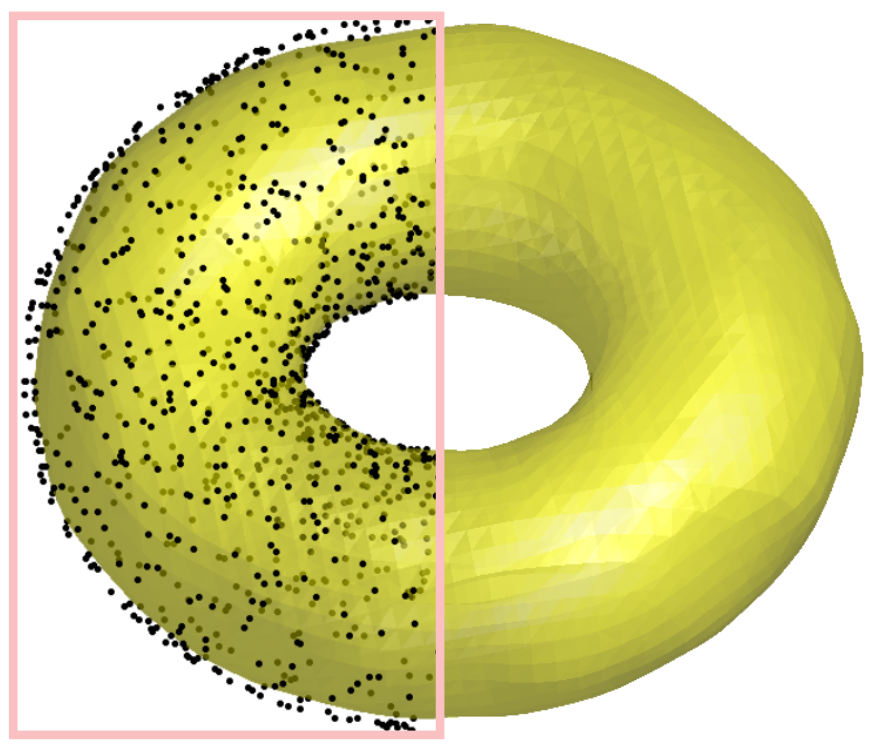}&
  \includegraphics[scale=0.35]{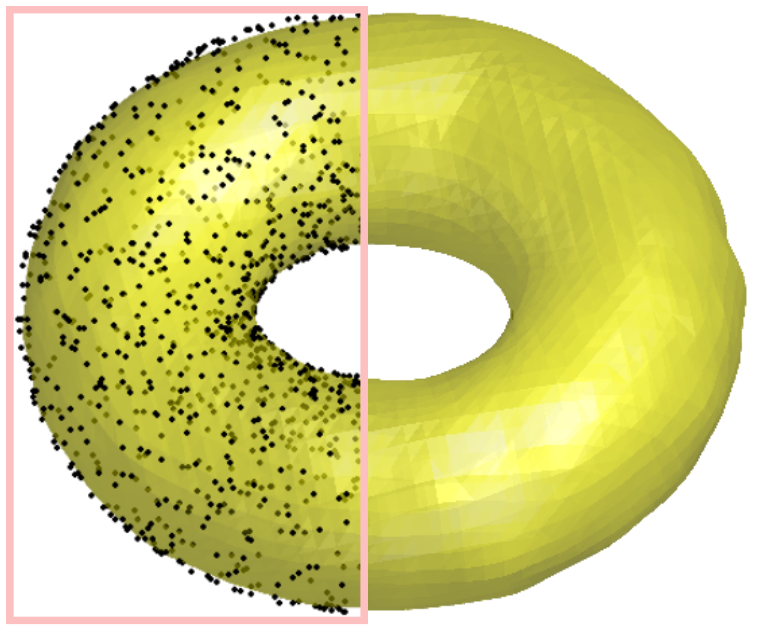}
   \end{tabular}
  \caption{The first row shows ALM and SIM applied to the 3D jar point cloud in Figure~\ref{fig::pcd}~(b). (a) The result of ALM with $r=1.3, \varepsilon=0.5, \eta=0.6$. (b) The result of SIM. The second row shows the methods applied to the 3D torus point cloud in Figure~\ref{fig::pcd}~(c). (c) The result of ALM with $r=1.3, \varepsilon=0.5, \eta=0.6$. (d) The result of SIM. Both methods are compatible and shows good results.}\label{fig.3d}
\end{figure}
The second set of experiments reconstruct surface from the point cloud $\D$ in 3D:  a jar in Figure~\ref{fig::pcd}~(b) and a torus in Figure~\ref{fig::pcd}~(c) within $\Omega=[0,50]^3$.   In Figure \ref{fig.3d}, we show  the reconstructed surfaces using SIM and ALM.  A portion of the given point cloud is superposed for validation in each case. Both methods successfully capture the overall shapes and non-convex features of the jar, as well as the torus.  There are  only slight differences in the reconstruction between using SIM with $p=2$ and using ALM with $p=1$.

Table \ref{tab.3d} shows the efficiency of SIM and ALM compared to the explicit method in \cite{ZOM*98} for the experiments in Figure \ref{fig.3d}.  Thanks to the semi-implicit scheme, $\Delta t$ can be large and we used $500$ in SIM;  in the explicit method, we are forced to use much smaller time step $\Delta t=20$ to maintain the stability.
The improvement of run-time in ALM is carefully  controlled by the parameters $r$, $\varepsilon$ and $\eta$. We choose $r=1.3$, $\varepsilon=0.5$ and $\eta=0.6$ for both cases. Both SIM and ALM efficiently provide accurate reconstruction.
\begin{table}[h]
  \centering
  \begin{tabular}{|c|c|c|c|}
    \hline
    % after \\: \hline or \cline{col1-col2} \cline{col3-col4} ...
    Object & ALM & SIM & \cite{ZOM*98} \\\hline
    Jar & $29.69$ & $29.42$ & $74.44$ \\\hline
    Torus & $47.32$ & $33.58$ & $114.20$ \\
    \hline
  \end{tabular}
  \caption{CPU time (s) of SIM and ALM compared to the explicit method in \cite{ZOM*98} for the  point cloud data sets of Figure~\ref{fig.3d}.  Both SIM and ALM show fast convergence. }\label{tab.3d}
\end{table}

\begin{figure}
\centering
\includegraphics[scale=0.45]{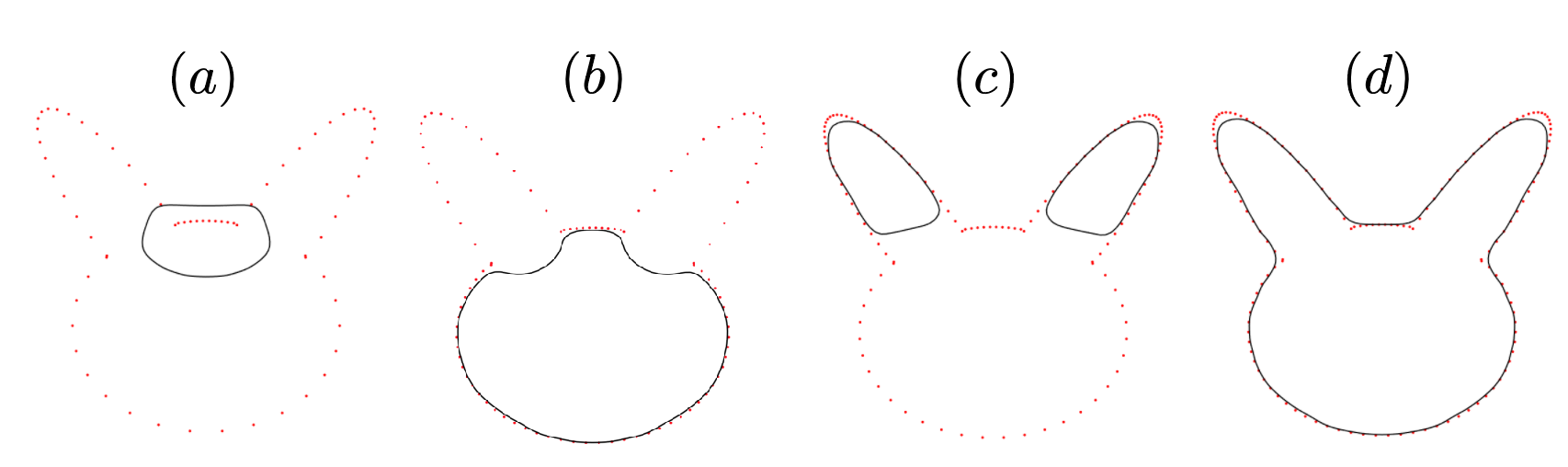}
\caption{The effect of the distance function for varying-density point clouds: the face with $n_1$ points, the head with $n_2$ points, and each ear with $n_3$ points.
(a) the given point cloud is with $(n_1, n_2, n_3)=(20,10,20)$, and shows the $0$-level-set of $\phi^n$ at $15^{\text{th}}$ iteration,  (b) $(n_1, n_2, n_3)=(50,10,20)$, and shows $18^{\text{th}}$ iteration, and (c) $(n_1, n_2, n_3)=(20,10,40)$, and shows $20^{\text{th}}$ iteration.  These three curves  eventually degenerate to a point. (d) is with $(n_1, n_2, n_3)=(50,10,40)$ and shows the converged solution.
The potential energy~(\ref{eq.energy0}) is mainly driven by the distance function $d$, which affects the level-set evolution.}\label{fig.pointdist}
\end{figure}

The third set of examples show the effect of the distance function $d$.  Notice that the weighted minimal surface energy~(\ref{eq.energy0}) is mainly driven by the distance function $d$, that is, the given point cloud $\D$ determines the landscape of $d$, which affects the behavior of the level-set during the evolution. Figure~\ref{fig.pointdist} shows the evolution using ALM, applied to different subsets of point clouds sampled from the same bunny face shape.  The densities of the point cloud vary
%are varied
for the three different regions: the face with $n_1$ points, the head with $n_2$ points, and each ear with $n_3$ points. Figure~\ref{fig.pointdist} (a) shows the given point cloud for $(n_1, n_2, n_3)=(20,10,20)$, with the $0$-level-set of $\phi^n$ at $15^{\text{th}}$ iteration,  (b) for $(n_1, n_2, n_3)=(50,10,20)$, at $18^{\text{th}}$ iteration, and(c) for $(n_1, n_2, n_3)=(20,10,40)$, at $20^{\text{th}}$ iteration.  These three curves  eventually degenerate to a point. (d) for $(n_1, n_2, n_3)=(50,10,40)$ and shows the converged solution.
In (a)--(c), denser parts of the point cloud attract the curve with stronger forces,  the sparser parts of the point cloud fail to lock the curve.  Then, the energy model~(\ref{eq.energy1}) drives curves to have short lengths, i.e., the level set tends to shrink.  In (d), with a more balanced distribution of points, the curve converges to correct shape.  %These examples are exadrrated to show the effect

\begin{figure}
  \centering
  \begin{tabular}{ccc}
   (a) & (b) & (c) \\
 \includegraphics[width=1.75 in]{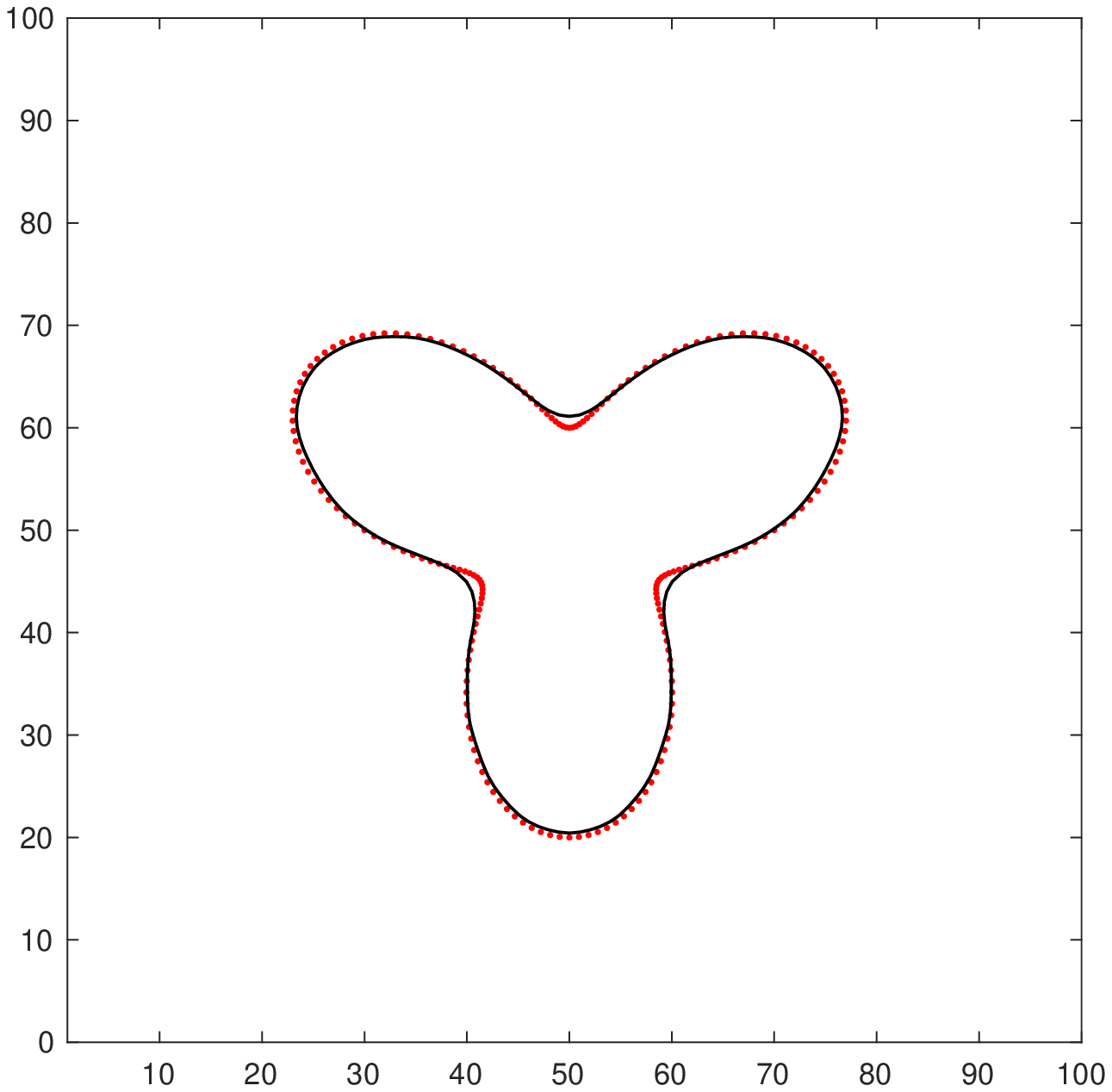} &
  \includegraphics[width=1.75in]{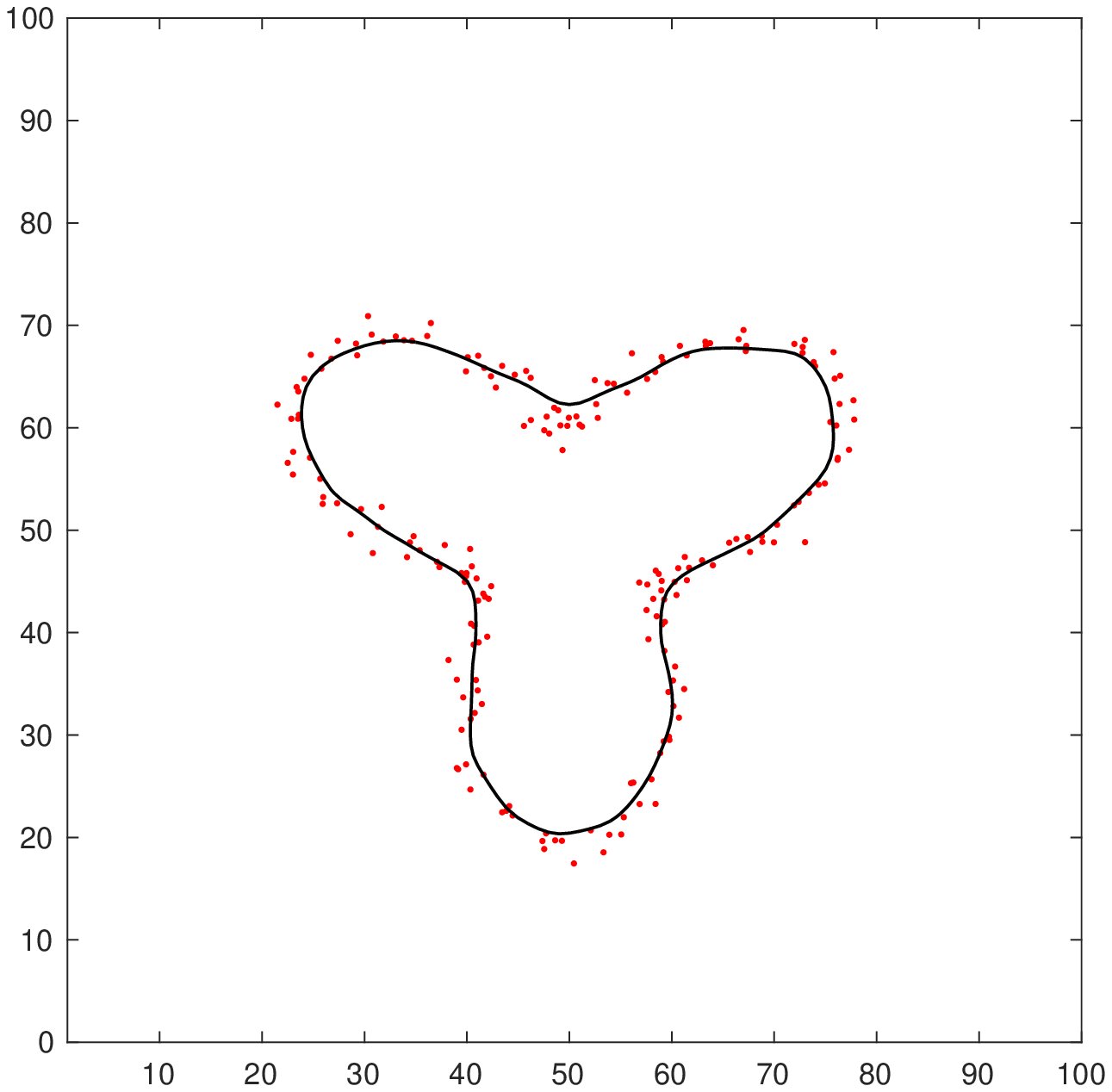} &
  \includegraphics[width=1.75 in]{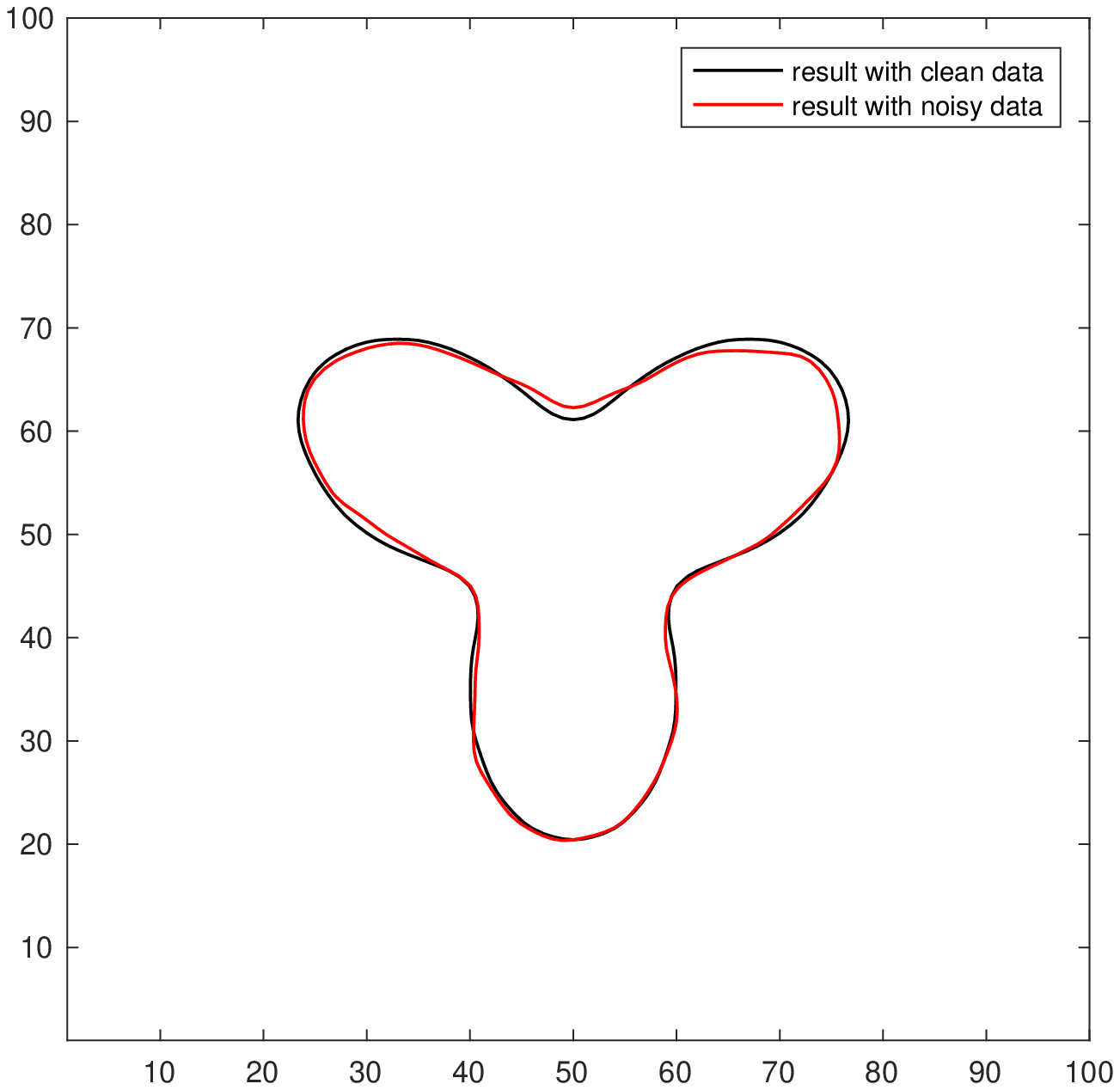}\\
  (d) & (e) & (f) \\
 \includegraphics[width=1.75 in]{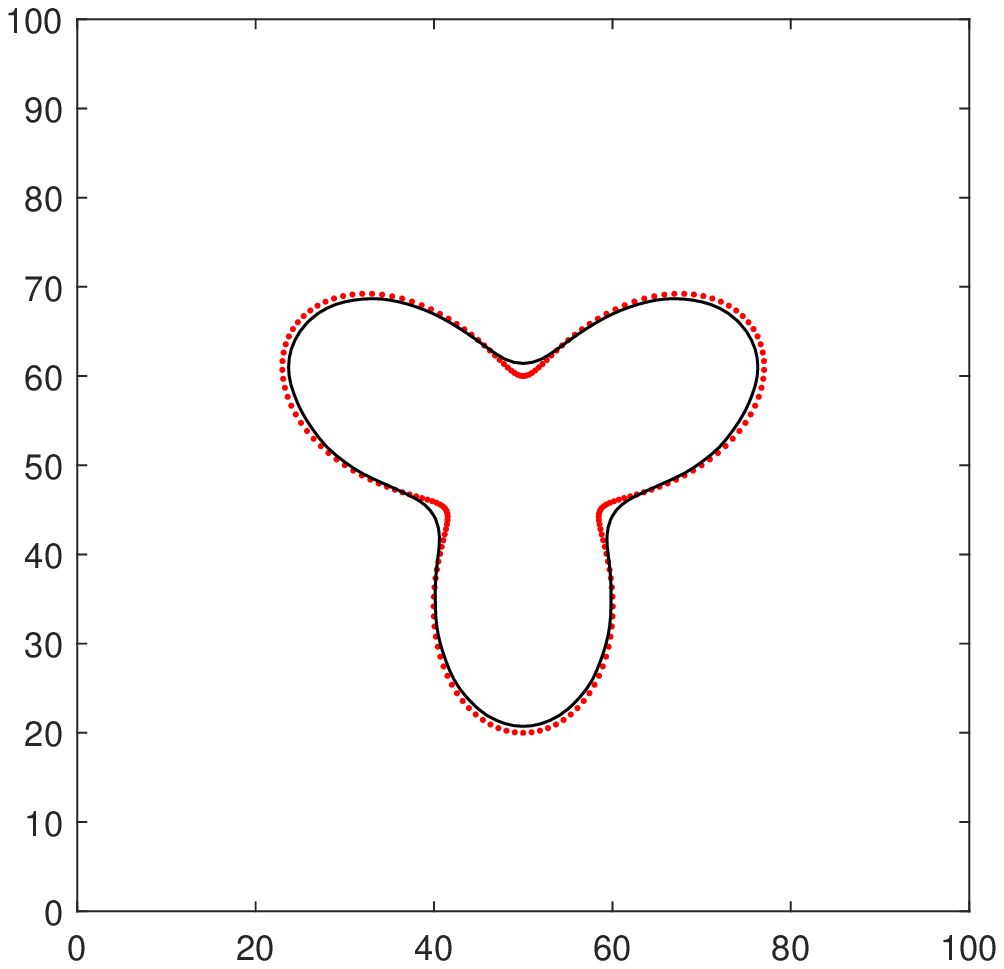} &
  \includegraphics[width=1.75 in]{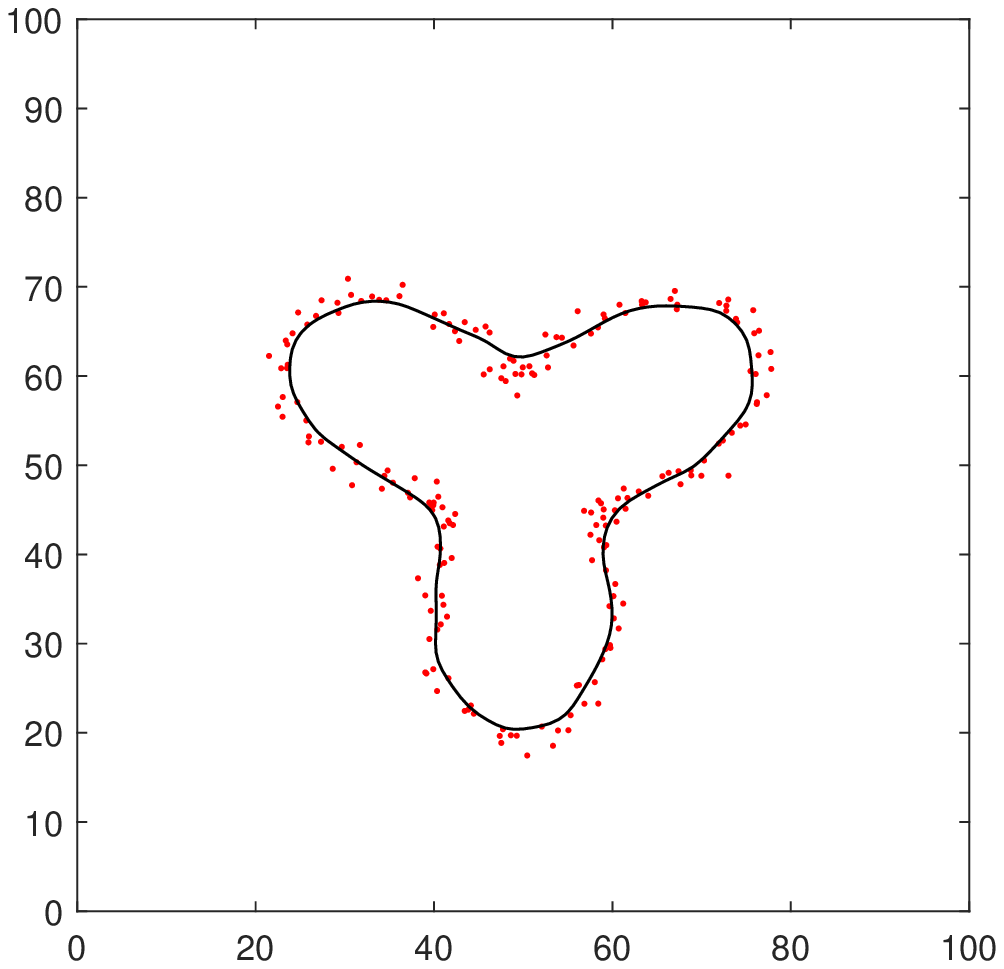} &
  \includegraphics[width=1.75 in]{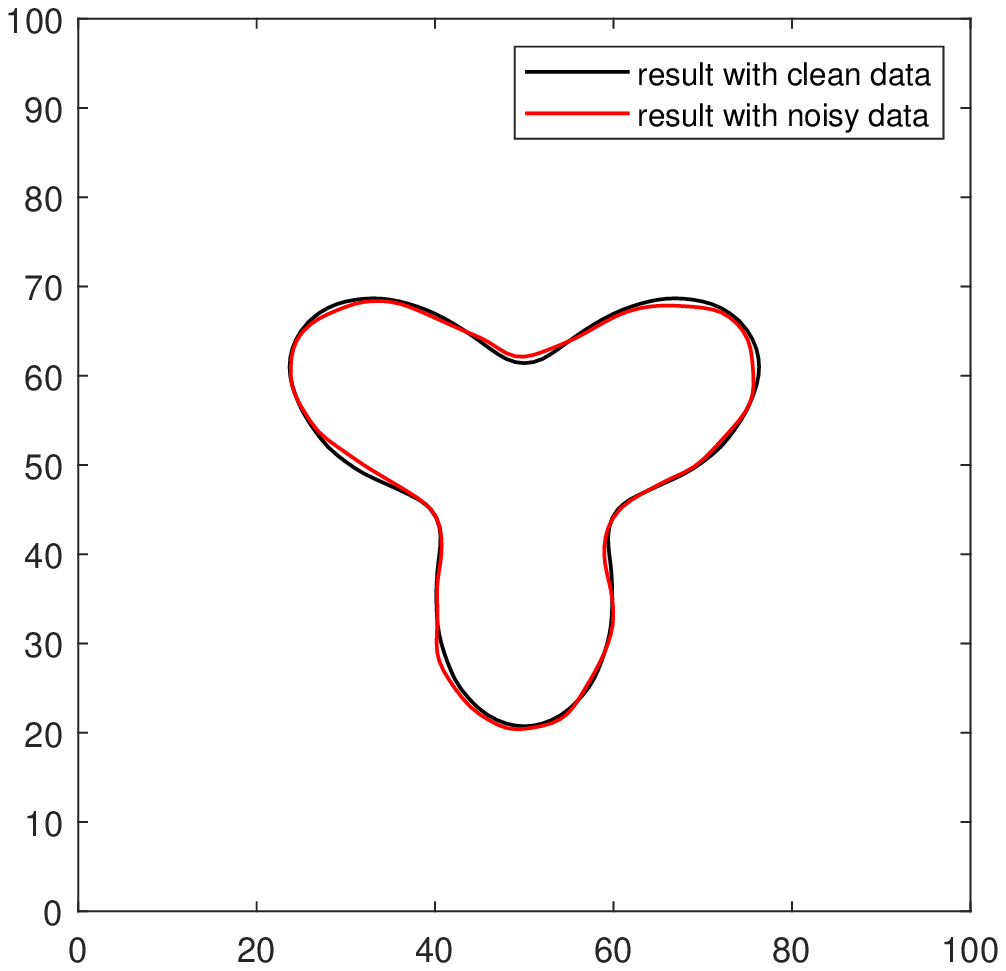}
  \end{tabular}
  \caption{The influence of noise on reconstructing three-fold circle with 200 points: (a)-(c) ALM and (d)-(f) SIM.  The first column shows the reconstructed curves from clean data, and the second column the reconstructions  from noisy data.  The third column shows the comparison between the two reconstructed curves in first two columns. }\label{fig.noise}
\end{figure}

The fourth set of examples demonstrate the robustness of ALM and SIM against noise.  Figure \ref{fig.noise} shows the reconstructed curves from clean and noisy data: (a)-(c) are results of ALM, and  (d)-(f)  are results of SIM.  (a) and (d) in the first column show results obtained from the clean data, which has 200 points sampled from a three-fold circle.  Gaussian noise with standard deviation $1$ is added to both $x$ and $y$ coordinates to generate noisy point cloud in the second column, (b) and (e).  To show the differences, the third column superposes both results reconstructed from clean and noisy point clouds. Both ALM and SIM  provide compatible results. For the noisy data, although the reconstructed curves show some oscillation, they are very close to the solutions using the clean data, respectively.
\subsection{Choice of Parameters for ALM and its Effects}\label{SEC::ALMbehavior} %\label{SEC::ParEffect}

The proposed ALM  has one parameter $r$, and the model (\ref{eq.energy1}) uses the delta function, where the smoothness parameter $\varepsilon$ is added to stabilize the computation.  Both parameters have straightforward effects on the level-set evolution from (\ref{eq.sub1_PDE}).  For example, consider a set of points within a thin-band around the 0-level-set of $\phi^n$, denoted by $B_{\varepsilon}=\{\x\mid -2\varepsilon/\sqrt{3}<\phi^n(\x)<2\varepsilon/\sqrt{3}\}$. By the continuity of $\phi^n$, there exist  $\mathbf{y}$ and $\mathbf{z}\in B_{\varepsilon}$ such that $\phi^n(\mathbf{y})=-\varepsilon/\sqrt{3}$ and $\phi^n(\mathbf{z})=\varepsilon/\sqrt{3}$; these values are the minimum and maximum of the function $h(x)=\frac{2\varepsilon x}{\pi(\varepsilon^2+x^2)^2}$, respectively. At these points, (\ref{eq.sub1_PDE}) takes the following forms:
\begin{align}
\Delta\phi^{n+1}=\begin{cases}
-9 \;d \; |\p^n|/(8\sqrt{3}\pi \;\varepsilon^2 r)+\nabla\cdot(\p^n+\frac{\bflambda^n}{r})\;\text{at~}\mathbf{y}.\\
9\; d\; |\p^n|/(8\sqrt{3}\pi \; \varepsilon^2 r)+\nabla\cdot(\p^n+\frac{\bflambda^n}{r})\;\;\;\;\text{at~}\mathbf{z}.
\end{cases}\label{eq::influence}
\end{align}
The first terms in the right hand side of (\ref{eq::influence}) show that with a smaller $\varepsilon$, there are less number of points in $B_{\varepsilon}$, but influence from $d$ is stronger.  With a larger $\varepsilon$, $d$ affects more number of points in $B_{\varepsilon}$, but the influence becomes weaker.  Varying $r$ also modifies the effect of $d$, while the size of $B_{\varepsilon}$ is not changed.

\begin{figure}
\centering
\includegraphics[scale=0.45]{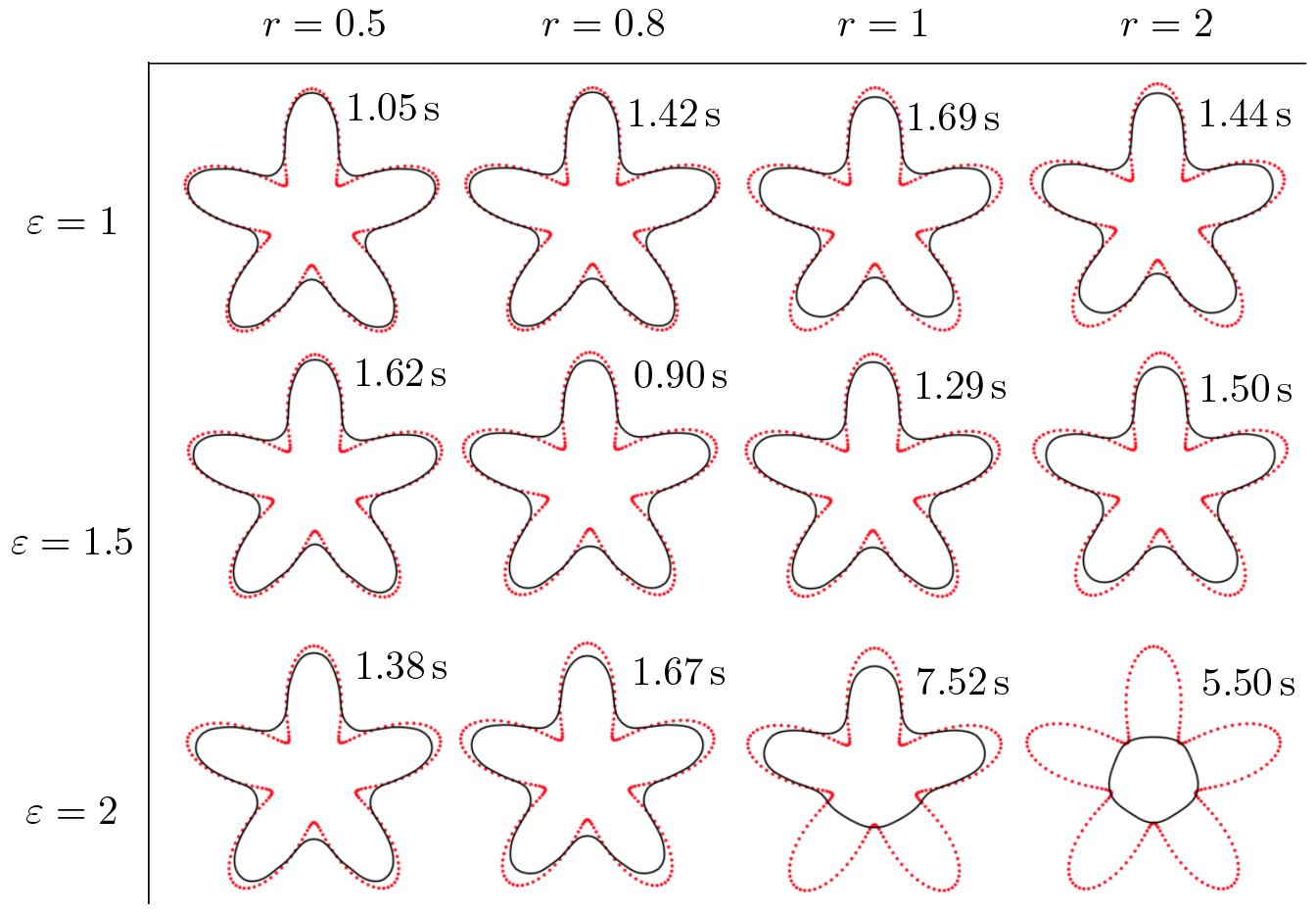}
\caption{Results by ALM with different $r$ and $\varepsilon$. For each column, from top to bottom, $\varepsilon=1,1.5,2$; and for each row, from left to right, $r=0.5,0.8,1,2$. Increasing $\varepsilon$ renders the curve less sharp and more convex.  Increasing $r$ induces a stronger diffusion effect on $\phi^n$.}\label{fig.fivecircle}
\end{figure}

We also find that $\varepsilon$ interacts with $r$ and effectively modifies the shape of the level-set.
Figure~\ref{fig.fivecircle} shows the results for ALM using different combinations of $r$ and $\varepsilon$, on the five-fold circle point cloud in Figure~\ref{fig::pcd} (a). For a fixed $r$, increasing $\varepsilon$ makes the approximated delta function smoother; consequently, narrow and elongated shapes are omitted, and the reconstructed curve becomes more convex.  For a fixed $\varepsilon$, larger $r$  loses more details, as discussed in Subsection \ref{SEC::conn}.   The speed of convergence varies for different combinations of $r$ and $\varepsilon$.  When  the choices are reasonable,  the algorithm converges fast within 2 seconds.  When both $r$ and $\varepsilon$ are large, results are not as good, and the convergences are slow.

Another observation comes from (\ref{eq.sub2_sol}). For any point $\x$ and $n\geq 0$, if the following value:
\begin{align*}
Q^n(\varepsilon,r) :=\phi^{n}\pi|r\nabla\phi^{n}-\bflambda^{n-1}|\varepsilon^2-d\varepsilon+(\phi^{n})^3\pi|r\nabla\phi^{n}-\bflambda^{n-1}|
\end{align*}
is positive, then  $\mathbf{p}^{n}(\x)=0$, and $d$ has no direct effect on (\ref{eq.sub1_PDE1}) at $\x$ in the next iteration.  Regarding $Q^n(\varepsilon,r)$ as a quadratic polynomial in terms of $\varepsilon$ parameterized by $r$,  the sign of $Q^n(\varepsilon,r)$ depends on the sign of $\phi^n(\x)$ and the sign of its discriminant computed via:
\begin{align*}
\text{Disc}\,Q^n=d^2-4(\phi^{n})^4\pi^2|r\nabla\phi^{n}-\bflambda^{n-1}|^2.
\end{align*}
%whose sign determines the number of roots of $Q^n$ for every fixed $r$.
The sign of $\phi^n(\x)$ is related to the position of $\x$ relative to the 0-level-set.  The sign of $\text{Disc}\,Q^n$ is determined by comparing the length of a vector difference $r\nabla\phi^{n}-\bflambda^{n-1}$ with the quantity $d/(4(\phi^n)^2\pi)$. By the projection theorem,  $|r\nabla\phi^{n}-\bflambda^{n-1}|^2$ is bounded below by $\alpha^n:=|\bflambda^{n-1}|^2-|\text{Proj}_{\nabla\phi^{n}}\bflambda^{n-1}|^2=|\bflambda^{n-1}|^2-|\bflambda^{n-1}\cdot\nabla\phi^n|^2/(|\bflambda^{n-1}|^2|\nabla\phi^n|^2)$, i.e., the squared residual of orthogonal projection of $\bflambda^{n-1}$ onto $\nabla\phi^n$; therefore, we can decide the sign of $\text{Disc}\,Q^n$ using $r$ via the following cases:
\begin{enumerate}
\item When $\frac{d^2}{4(\phi^{n})^4\pi^2}< \alpha^n$,
%where  $\alpha^n=|\bflambda^{n-1}|^2-|\text{Proj}_{\nabla\phi^{n}}\bflambda^{n-1}|^2$ with $\text{Proj}_\mathbf{v}\mathbf{u}=\frac{\mathbf{u}\cdot\mathbf{v}}{|\mathbf{v}|}\frac{\mathbf{v}}{|\mathbf{v}|}$ as the orthogonal projection of $\mathbf{u}$ onto $\mathbf{v}$,
for any $r>0$, $\text{Disc}\,Q^n<0$.
\item When $\frac{d^2}{4(\phi^{n})^4\pi^2}\geq \alpha^n$:
\begin{enumerate}
\item if $r>r_U^n$ or $r<r_L^n$, then $\text{Disc}\,Q^n<0$;
\item if $\max\{0,r_L^n\}\leq r\leq r_U^n$, then $\text{Disc}\,Q^n\geq0$.
\end{enumerate}
Here,
\begin{equation*}%\label{eq::rcondi}
r^n_U=\frac{|\text{Proj}_{\nabla\phi^{n}}\bflambda^{n-1}|+\sqrt{\frac{d^2}{4(\phi^{n})^4\pi^2}-\alpha^n}}{|\nabla\phi^{n}|} \;\; \text{ and } \;\;  r^n_L=\frac{|\text{Proj}_{\nabla\phi^{n}}\bflambda^{n-1}|-\sqrt{\frac{d^2}{4(\phi^{n})^4\pi^2}-\alpha^n}}{|\nabla\phi^{n}|}.
\end{equation*}
\end{enumerate}
%This is computed by considering (\ref{disc}) as a difference between two vectors and finding the expression in terms of $r$.
%
%\begin{enumerate}
%\item When $\text{Disc}\,Q^n<0$,
%\begin{enumerate}
%\item {$r$ can be any positive real number, if $\frac{d^2}{4(\phi^{n})^4\pi^2}< \alpha^n$. Here  $\alpha^n=|\bflambda^{n-1}|^2-|\text{Proj}_{\nabla\phi^{n}}\bflambda^{n-1}|^2$ with $\text{Proj}_\mathbf{v}\mathbf{u}=\frac{\mathbf{u}\cdot\mathbf{v}}{|\mathbf{v}|}\frac{\mathbf{v}}{|\mathbf{v}|}$ as the orthogonal projection of $\mathbf{u}$ onto $\mathbf{v}$;}
%\item{When $\frac{d^2}{4(\phi^{n})^4\pi^2}> \alpha^n$,  $r>r^n_U$ or $ r<r^n_L$, where
%if $r_L^n$ is positive.  }
%\end{enumerate}
%\item When $\text{Disc}\,Q^n\geq 0$,
%\begin{enumerate}
%	\item{there exists no such $r$, if $\frac{d^2}{4(\phi^{n})^4\pi^2}< \alpha^n$, and  }
%	\item{ when $\frac{d^2}{4(\phi^{n})^4\pi^2}> \alpha^n$, $r$ satisfies
%$	\max\{0,r_L^n\}\leq r\leq r_U^n.$}
%\end{enumerate}
%\end{enumerate}
When $\phi^n(\x)>0$, $Q^n$ concaves upwards and $Q^n(0,r)\geq0$ for any $r$. If $\text{Disc}\,Q^n<0$, $Q^n$ is positive for all $\varepsilon$ and $d$ has no effect on level set evolution.   If $\text{Disc}\,Q^n\geq0$,  $Q^n$ is  positive for $\varepsilon$ outside the interval bounded by two roots of $Q^n$, i.e.,
\begin{align*}
0<\varepsilon<\frac{d-\sqrt{\text{Disc}\,Q^n}}{2\phi^n\pi|r\nabla\phi^n-\bflambda^{n-1}|} \;\; \text{ or } \;\;\varepsilon>\frac{d+\sqrt{\text{Disc}\,Q^n}}{2\phi^n\pi|r\nabla\phi^n-\bflambda^{n-1}|}\;.
\end{align*}
When $\phi^n(\x)<0$, $Q^n$ concaves downwards, and $Q^n(0,r)\leq 0$ for any $r$. In this case, $Q^n$ is never positive: either $\text{Disc}\,Q^n<0$, i.e., no roots, or $\text{Disc}\,Q^n\geq 0$ but both roots are negative.

Notice that the bounds, $r_L^n$ and $r_U^n$, are closely related to the ratio $d/(\phi^n)^2$, which contributes to the adaptive behavior of ALM.    For example, for a point $\x$ where $\phi^n(\x)>0$, when $|\phi^n(\x)|$ is close to $0$ but $d(\x)\gg 0$, $r_L^n<0$ and $r_U^n$ becomes extremely large; thus, for a moderate value of $r$, $d$ has strong influence on the evolution of the level-set near $\x$ and swiftly moves the curve towards the point cloud.  For a point $\x$ which is close to both $\D$ and $\{\phi^{n}=0\}$, the level-set evolution becomes more stringent about the minimization of the energy~(\ref{eq.energy2}).

\begin{figure}
	\begin{tabular}{cc}
		(a)&(b)\\	
		\includegraphics[width=2.65in]{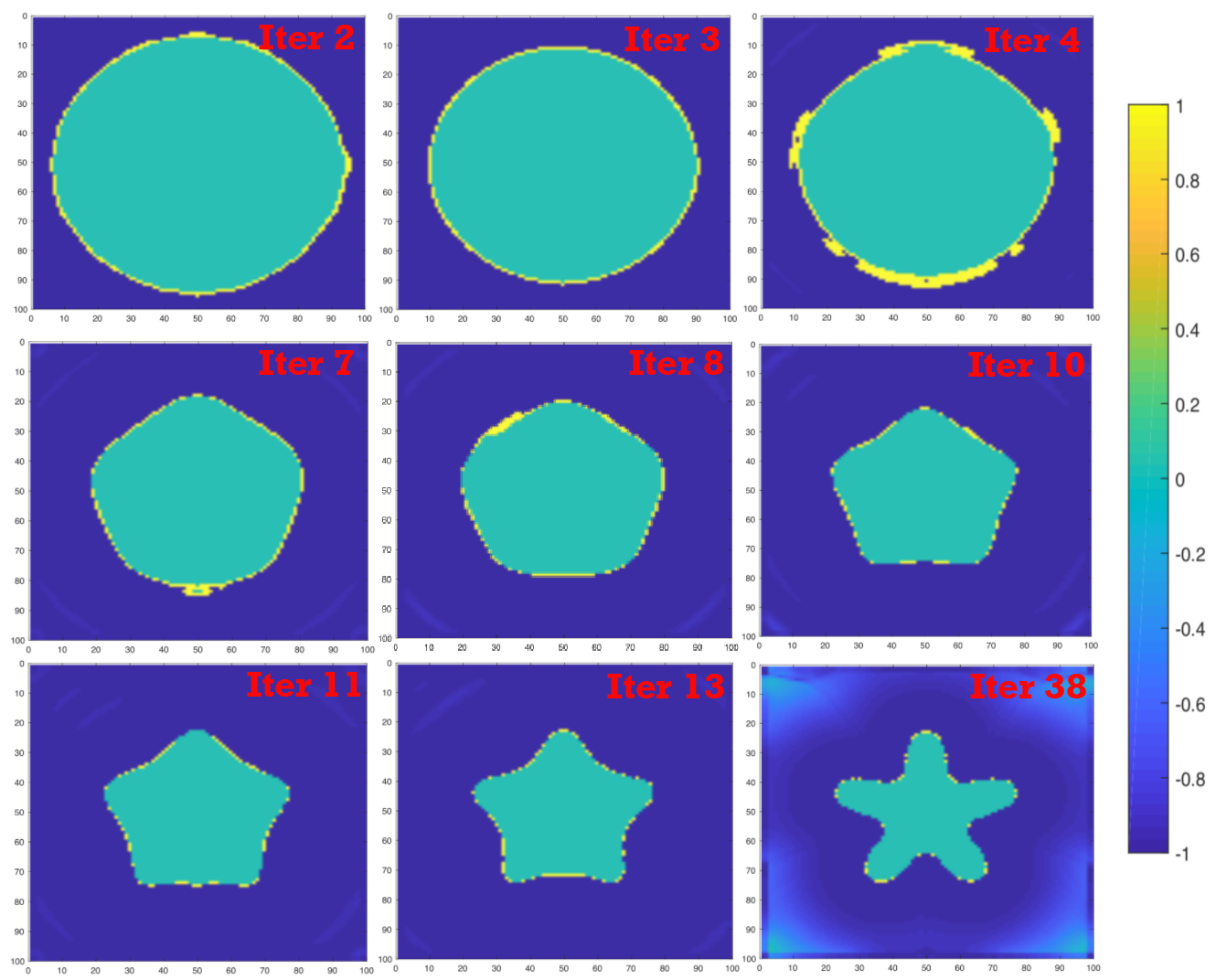} &	
		\includegraphics[width=2.65in]{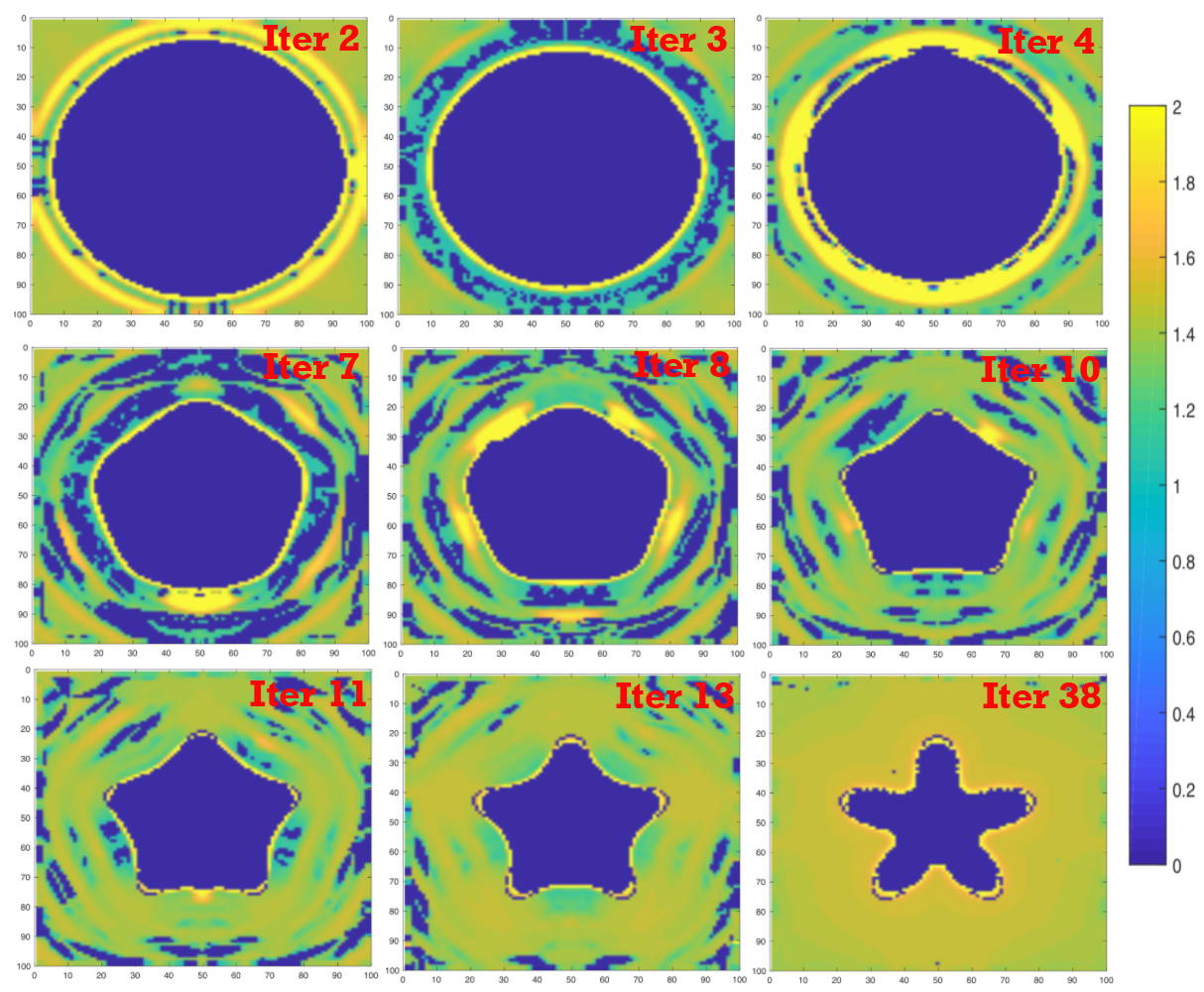} \\
		(c)&(d)\\
		\includegraphics[width=2.65in]{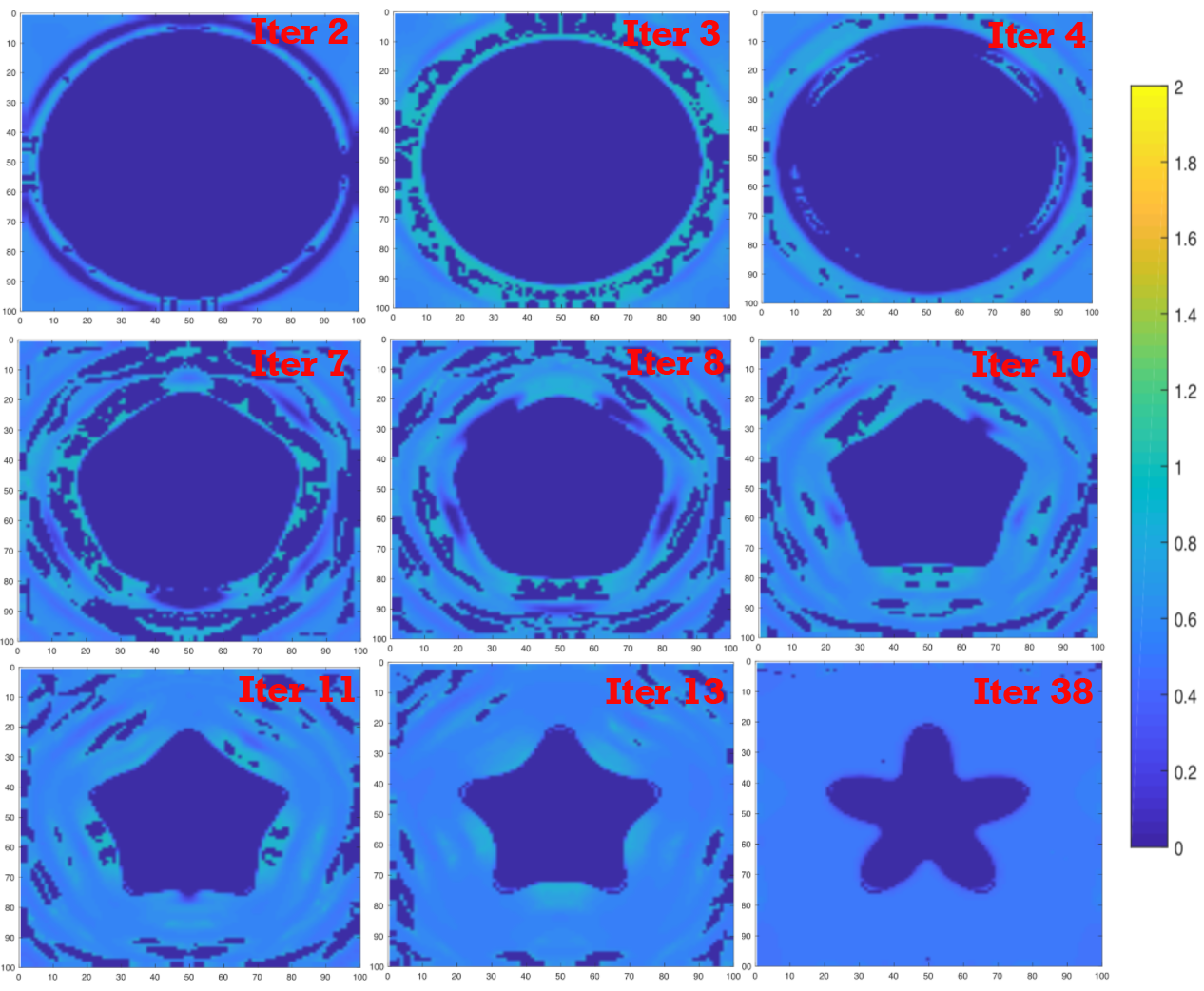} &
		\includegraphics[width=2.45in]{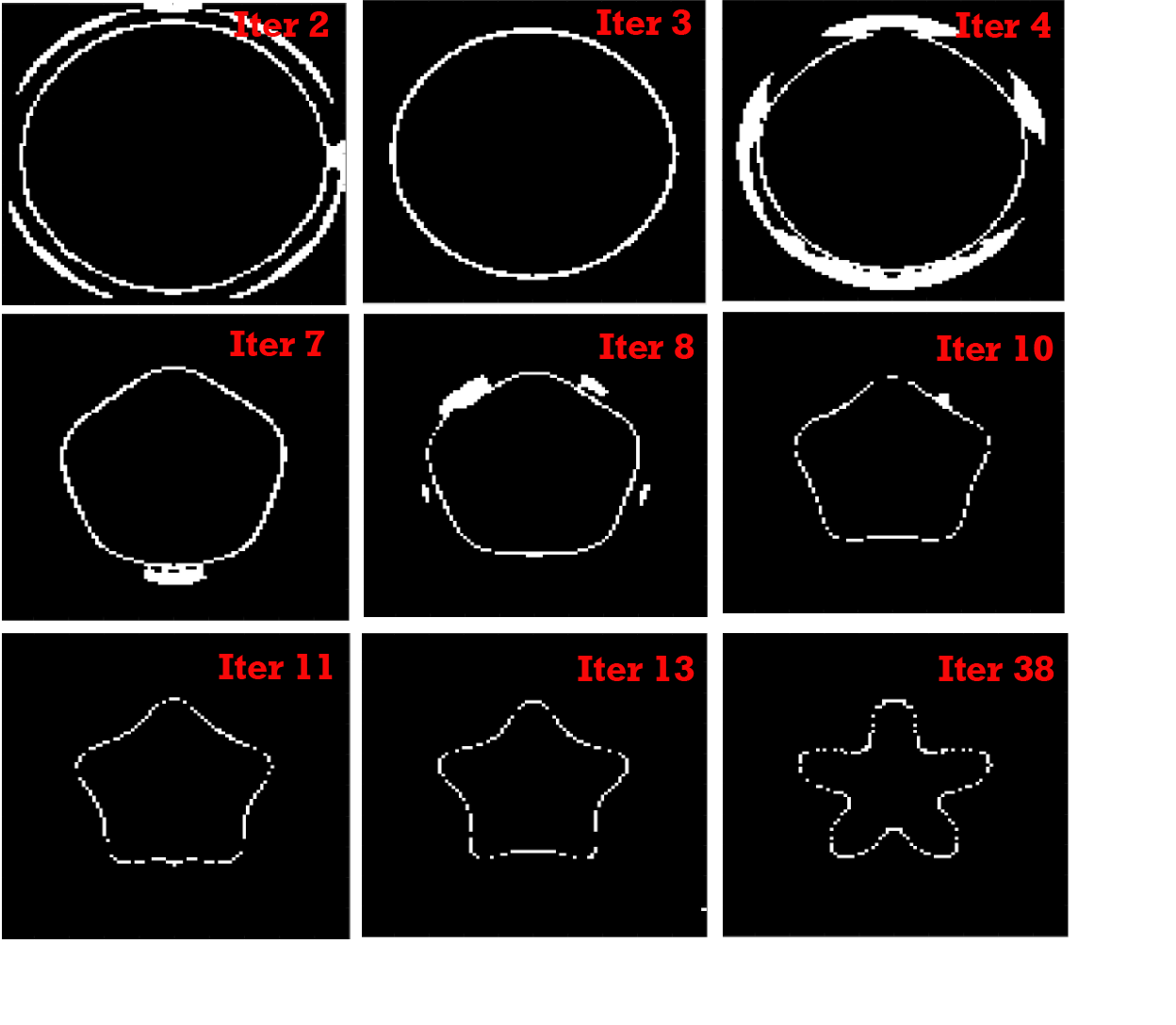}
	\end{tabular}
	\caption{(a) $\text{Disc}\,Q^n$, (b) $r_U^n$, (c) $r_L^n$ at certain iterations. (d) The region (in white) where $d$ explicitly guides the level-set evolution by ALM. The distance function $d$ refines the local structures and it is only active near $\{\phi^n=0\}$. This partially explains the efficiency of ALM. }\label{fig::rbound}
\end{figure}
Figure~\ref{fig::rbound} illustrates this effect, for the five-fold circle point cloud in Figure~\ref{fig::pcd}~(a) with $r=2$ and $\varepsilon=1$. Figure~\ref{fig::rbound} shows (a) $\text{Disc}\,Q^n$, (b) $r^n_U$, (c) $r^n_L$ , and (d)  the region where $d$ effects the level set evolution.  The figures are  for iterations $n=2,3,4,7,8,10,11,13$ and $38$ (converged).   The region inside $\{\phi^n=0\}$ always experiences the influence of $d$, as described above.  Figure~\ref{fig::rbound}~(a) shows that the region outside $\{\phi^n=0\}$ is mostly blue indicating $\text{Disc}\,Q^n<0$; hence, for almost every point outside the 0-level-set, as long as $r^n_L\leq r\leq r_U^n$, the landscape of $d$ has strong effects on the evolution. In (b) and (c), observe that high values of $r^n_U$ only concentrate near the 0-level-set while $r_L^n$ remains relatively small in the whole domain;  thus, %, according to (\ref{eq::rcondi}),
the influence of $d$ is strong near $\{\phi^n=0\}$. (d) displays the white regions where $d$ explicitly guides the level-set evolution and the black regions where $d$ has no direct effect. These results show that, although ALM evolves the level-set globally, it ignores the effects of $d$ when evolving the regions far away from the level-sets; and it utilizes $d$ to refine the local structures for the regions of the level-sets close to $\D$. %This partially explains the efficiency of ALM.

\section{Conclusion}
\label{SEC:concl}

We propose two fast algorithms, SIM and ALM, to reconstruct $m$-dimensional manifold from unstructured point clouds by minimizing the weighted minimum surface energy (\ref{eq.energy1}). SIM improves the computational efficiency by relaxing the constraint on the time-step using a semi-implicit scheme. ALM follows an augmented Lagrangian approach and solves the problem by an ADMM-type algorithm. Numerical experiments show that the proposed algorithms are superior at the computational speed, and both of them produce accurate results. Theoretically, we demonstrate the delicate interaction among parameters involved in ALM, and show the connections between SIM and ALM.  This explains the behaviors of ALM from the perspective of SIM.

%\bibliographystyle{abbrv}
%\bibliography{cite_fastSurfR}

\end{document}